\documentclass[11pt,a4paper,leqno,twoside]{amsart}
\usepackage{latexsym,amssymb,amsmath}
\input xy
\xyoption{all}

\def\CC{\mathbb C}

\def\RR{\mathbb R}
\def\HH{\mathbb H}
\def\AA{{\mathbb A}}

\def\OO{\mathbb O}

\def\11{\mathbf 1}
\def\PP{\mathbb P}

\def\e1{\varepsilon_1}
\def\e2{\varepsilon_2}
\def\e3{\varepsilon_3}

\def\P2{{\PP}^2}

\def\00{\underline{0}}
\def\J0{{\cal J}_3(\underline{0})}

\def\PJ0{\PP({\cal J}_3(\underline{0}))}

\def\e{\varepsilon}

\def\AP2{{\AA\PP}^2}
\def\RP2{{\RR\PP}^2}
\def\CP2{{\CC\PP}^2}
\def\HP2{{\HH\PP}^2}
\def\OP2{{\OO\PP}^2}

\newtheorem{theo}{Theorem}[section]
\newtheorem{coro}[theo]{Corollary}
\newtheorem{lemm}[theo]{Lemma}
\newtheorem{prop}[theo]{Proposition}

\theoremstyle{remark}
\newtheorem{rema}[theo]{Remark}

\begin{document}
\title[Quasifinite fields of Diophantine dimension]{Quasifinite
fields of prescribed characteristic and Diophantine dimension}

\author{Ivan D. Chipchakov}
\address{Institute of Mathematics and Informatics\\Bulgarian Academy
of Sciences\\Acad. G. Bonchev Str., Block 8\\1113 Sofia, Bulgaria:
E-mail: chipchak@math.bas.bg}

\author{Boyan B. Paunov}
\address{Faculty of Mathematics and Informatics\\Sofia
University\\ James Bourchier Blvd. 5\\1164 Sofia, Bulgaria: E-mail
address: bpaunov@fmi.uni-sofia.bg}

\begin{abstract}
Let $\mathbb{P}$ be the set of prime numbers, $\overline {
\mathbb{P}}$ the union $\mathbb{P} \cup \{0\}$, and for any field
$E$, let char$(E)$ be its characteristic, ddim$(E)$ the
Diophantine dimension of $E$, $\mathcal{G}_{E}$ the absolute
Galois group of $E$, and cd$(\mathcal{G}_{E})$ the Galois
cohomological dimension $\mathcal{G}_{E}$. The research presented
in this paper is motivated by the open problem of whether
cd$(\mathcal{G}_{E}) \le {\rm ddim}(E)$. It proves the existence
of quasifinite fields $\Phi _{q}\colon q \in \mathbb{P}$, with
ddim$(\Phi _{q})$ infinity and char$(\Phi _{q}) = q$, for each
$q$. It shows that for any integer $m > 0$ and $q \in \overline
{\mathbb{P}}$, there is a quasifinite field $\Phi _{m,q}$ such
that char$(\Phi _{m,q}) = q$ and ddim$(\Phi _{m,q}) = m$. This is
used for proving that for any $q \in \overline {\mathbb{P}}$ and
each pair $k$, $\ell \in (\mathbb{N} \cup \{0, \infty \})$
satisfying $k \le \ell $, there exists a field $E _{k, \ell ; q}$
with char$(E _{k, \ell ; q}) = q$, ddim$(E _{k, \ell ; q}) = \ell
$ and cd$(\mathcal{G}_{E_{k, \ell ; q}}) = k$. Finally, we show
that the field $E _{k, \ell ; q}$ can be chosen to be perfect
unless $k = 0 \neq \ell $.
\end{abstract}

\maketitle

\par
\medskip
\section{\bf Introduction}

\medskip

Let $F$ be a field, $F _{\rm sep}$ its separable closure, Fe$(F)$
the set of finite extensions of $F$ in $F _{\rm sep}$, and
$\mathcal{G}_{F}$ the absolute Galois group of $F$, i.e. the
Galois group $\mathcal{G}(F _{\rm sep}/F)$. By definition, the
Brauer dimension Brd$(F)$ of $F$, introduced in \cite{ABGV}, is
equal to the least integer $n \ge 0$, for which the degree
deg$(D)$ of any finite-dimensional associative central division
$F$-algebra $D$ divides exp$(D) ^{n}$, where exp$(D)$ is the
exponent of $D$, i.e. the order of the equivalence class $[D]$ of
$D$ as an element of the Brauer group Br$(F)$; if no such $n$
exists, Brd$(F)$ is defined to be infinity. The Brauer
$p$-dimension Brd$_{p}(F)$ is defined analogously, for each prime
number $p$, by letting $D$ run across the class of central
division $E$-algebras of $p$-primary degrees. In view of the
primary tensor product decomposition theorem for central division
algebras over fields (cf. \cite{P}, Sect. 14.4), Brd$(F)$ equals
the supremum of Brd$_{p}(F)\colon p \in \mathbb{P}$, where
$\mathbb{P}$ is the set of prime numbers. As in \cite{Ch1}, the
supremum sup$\{{\rm Brd}_{p}(Y)\colon Y \in {\rm Fe}(F)\}$,
denoted by abrd$_{p}(F)$, is called an absolute Brauer
$p$-dimension of $F$. We say that $F$ is a field of dimension $\le
1$, in the sense of Serre (see \cite{S1}, Ch. II, 3.1), if Br$(Y
^{\prime }) = \{0\}$, i.e. Brd$(Y ^{\prime }) = 0$, for every
algebraic field extension $Y ^{\prime }/F$; in view of the
Albert-Hochschild theorem \cite[Ch. II, 2.2]{S1} (and well-known
general properties of maximal separable subextensions of algebraic
extensions, see \cite[Ch. V, Sects. 4 and 6]{L}), this holds if
and only if abrd$_{p}(F) = 0$, for every $p \in \mathbb{P}$.
\par
The cohomological dimension cd$(\mathcal{G}_{F})$ of
$\mathcal{G}_{F}$ (viewed as a profinite group, whence, a compact
totally disconnected group with respect to the Krull topology) is
defined to be the supremum of the cohomological $p$-dimensions
cd$_{p}(\mathcal{G}_{F})$, $p \in \mathbb{P}$. By Galois
cohomology (cf. \cite{S1}, Ch. I, Proposition~21),
cd$_{p}(\mathcal{G}_{F}) \le n$, for a given $p \in \mathbb{P}$
and an integer $n \ge 0$, if and only if the cohomology group $H
^{n+1}(\mathcal{G}_{F}, \mathbb{Z}/p\mathbb{Z})$ is trivial, where
the field $\mathbb{Z}/p\mathbb{Z}$ (of residue classes of integers
modulo $p$) is viewed as a discrete $\mathcal{G}_{F}$-module with
the trivial action of $\mathcal{G}_{F}$; also,
cd$_{p}(\mathcal{G}_{F}) = n$ if and only if $n$ is minimal with
this property. We set cd$_{p}(\mathcal{G}_{F}) = \infty $ if $H
^{n}(\mathcal{G}_{F}, \mathbb{Z}/p\mathbb{Z}) \neq \{0\}$, for all
$n \in \mathbb{N}$. When $p \neq {\rm char}(F)$, it follows that
cd$_{p}(\mathcal{G}_{F}) \le n < \infty $ if and only if the
cohomology group $H ^{n+1}(\mathcal{G}_{F(\mu (p))}, \mu (p))$ is
trivial, $F(\mu (p))$ being the extension of $F$ generated by the
set (in fact, a multiplicative group) $\mu (p)$ of $p$-th roots of
unity lying in $F _{\rm sep}$ (cf. \cite[Ch. I,
Proposition~14]{S1}, and \cite[Ch. VI, Sect. 3]{L}). Thus the
computation of cd$_{p}(\mathcal{G}_{F})$ reduces to the special
case in which $F$ contains a primitive $p$-th root of unity unless
$p = {\rm char}(F)$.
\par
It is well-known (cf. \cite[Ch. II, 3.1]{S1}) that dim$(F) \le 1$
if and only if Br$(F ^{\prime }) = \{0\}$, where $F ^{\prime }$
runs across Fe$(F)$. For example, dim$(F _{\rm sep}) \le 1$; this
follows from the Albert-Hochschild theorem if char$(F) \neq 0$.
When $F$ is perfect, we have dim$(F) \le 1$ if and only if
cd$(\mathcal{G}_{F}) \le 1$. For example, dim$(F) \le 1$ if $F$ is
a quasifinite field, i.e. a perfect field which admits a unique
extension in $F _{\rm sep}$ of degree $n$, for each $n \in
\mathbb{N}$. It is known that then $\mathcal{G}_{F}$ is isomorphic
to $\mathcal{G}_{\mathbb{F}}$, for any finite field $\mathbb{F}$.
The inequality cd$(\mathcal{G}_{F}) \le 1$ holds, since
$\mathcal{G}_{F}$ is isomorphic to the topological group product
$\prod _{p \in \mathbb{P}} \mathbb{Z} _{p}$, where $\mathbb{Z}
_{p}$ is the additive group of $p$-adic integers, for each $p$
(see \cite[Examples~4.1.2]{GiSz}).
\par
We say that the Diophantine dimension ddim$(F)$ of $F$ is finite
and equal to $m$, if $m$ is the least integer $\ge 0$, for which
$F$ is a field of type $C _{m}$; if no such $m$ exists, we set
ddim$(F) = \infty $. By type $C _{m}$ (or a $C _{m}$-field), for
an integer $m \ge 0$, we mean that every $F$-form (a homogeneous
nonzero polynomial with coefficients in $F$) $f$ of degree
deg$(f)$ in more than deg$(f) ^{m}$ variables has a nontrivial
zero over $F$. For example, $F$ is a $C _{0}$-field, i.e. ddim$(F)
= 0$ if and only if it is algebraically closed. The class of $C
_{m}$-fields is closed under taking algebraic extensions, and it
contains every extension of transcendence degree $m$ over any
algebraically closed field (cf. \cite{L1}). The question of
whether this class consists of fields of Brauer dimensions less
than $m$ is presently open; it is known, however, that these
fields have absolute Brauer $p$-dimensions less than $p ^{m-1}$,
for all $p \in \mathbb{P}$ (see \cite{Mat}). When $F$ is a $C
_{m}$-field and char$(F) =  q > 0$, it is easily verified that
$[F\colon F ^{q}] \le q ^{m}$, where $F ^{q} = \{\alpha ^{q}\colon
\alpha \in F\}$ (the $F$-form $\sum _{i=1} ^{q^{m'}} a _{i}X _{i}
^{q}$ does not possess a nontrivial zero over $F$, provided
$[F\colon F ^{q}] = q ^{m'}$ and the system $a _{i} \in F\colon i
= 1, \dots , q ^{m'}$, is a basis of $F$ over $F ^{q}$). Since, by
Galois cohomology (cf. \cite[Ch. I, 3.3]{S1}),
cd$(\mathcal{G}_{F}) = 0$ if and only if $F = F _{\rm sep}$, these
observations prove the following:
\par
\medskip\noindent
(1.1) In order that cd$(\mathcal{G}_{F}) = 0$ and ddim$(F) = m <
\infty $, it is sufficient that $F = F _{\rm sep}$, char$(F) = q >
0$, $[F\colon F ^{q}] = q ^{m}$, and $F$ is an extension of an
algebraically closed field $F _{0}$ of transcendence degree $r$.
We have cd$(\mathcal{G}_{F}) = 0$ and ddim$(F) = \infty $ if $F =
F _{\rm sep}$, char$(F) = q > 0$ and $[F\colon F ^{q}] = \infty $.
\par
\medskip
The research presented in this paper is motivated by the following
open question (arising from an observation made by Serre in
\cite[Ch. II, 4.5]{S1}):
\par
\medskip
{\bf Question 1.} Find whether cd$(\mathcal{G}_{F}) \le m$
whenever $F$ is a field of type $C _{m}$, for some $m \in
\mathbb{N}$. Equivalently, find whether cd$(\mathcal{G}_{F}) \le
{\rm ddim}(F)$.
\par
\medskip
Finite fields have type $C _{1}$, by the classical
Chevalley-Warning theorem, and $C _{1}$-fields have dimension $\le
1$ (see \cite[Theorem~2.6]{GiSz}, and \cite[Ch. II, 3.1 and
3.2]{S1}, respectively). Nontrivially, it follows from
Merkur'ev-Suslin's theorem (cf. \cite[Corollary~24.9]{Su}) that
cd$(\mathcal{G}_{E}) \le 2$, for every $C _{2}$-field $E$.
Question 1 is open if $m \ge 3$. It has been proved that
cd$_{p}(\mathcal{G}_{E_{m}}) < \infty $, $p \in \mathbb{P}$, for
every $C _{m}$-field $E _{m}$; however, in case $m \ge 3$, the
sequence $c _{p}(m)$, $p \in \mathbb{P}$, of best presently known
explicit upper bounds on cd$_{p}(\mathcal{G}_{E _{m}})$, found in
\cite{KM}, is unbounded, which does not rule out the possibility
that cd$(\mathcal{G}_{E _{m}}) = \infty $.
\par
\medskip
Note that if char$(F) = q > 0$ and $F$ is a $C _{m}$-field, then
the $q$-dimension dim$_{q}(F)$ introduced, for example, in
\cite{KK}, is at most equal to $m$. This theorem, established by
Arason and Baeza \cite{ArBa}, implies the answer to Question 1
will be affirmative if and only if ddim$(F) \ge {\rm CD}(F)$, for
every field $F$, where CD$(F)$ is the cohomological dimension of
$F$ (in the sense of \cite{IzAr}), defined as follows:
\par
\medskip
{\bf Definition 1.} {\it {\rm (i)} {\rm CD}$(F) = {\rm
cd}(\mathcal{G}_{F})$ if {\rm char}$(F) = 0$.
\par
{\rm (ii)} When {\rm char}$(F) = q > 0$, {\rm CD}$(F)$ equals the
supremum of {\rm dim}$_{q}(F)$ and {\rm
cd}$_{p}(\mathcal{G}_{F})\colon p \in \mathbb{P} \setminus \{q\}$;
{\rm dim}$_{q}(F)$ is defined to be the minimal integer $u$ for
which $[F\colon F ^{q}] \le q ^{u}$ and the Kato-Milne cohomology
group $H _{q} ^{u+1}(L)$ is trivial, for every finite extension
$L$ of $F$ (or to be infinity if such no $u$ exists).}
\par
\medskip
For the definition of the groups $H _{q} ^{n+1}(F)$, $n \ge 0$, in
characteristic $q$, and for more information about them, we refer
the reader to \cite{KK} and \cite{IzAr}. An interpretation of
these groups as the $q$-part of the Galois cohomology of $F$ has
been made in \cite{Izh}, and results on the minimal index $n$ for
which $H _{q} ^{n+1}(F) = \{0\}$ can be found, e.g., in
\cite{ArArBa}. Here we note that $H _{q} ^{2}(F)$ is isomorphic to
the maximal subgroup of Br$(F)$ of order dividing $q$ (cf.
\cite[Sect. 9.2]{GiSz} or \cite[pages~219-220]{Ka}). This implies
dim$_{q}(F) \le 1$ if and only if Br$(L)$ does not contain an
element of order $q$, for any finite extension $L/F$; in addition,
it follows that CD$(F) \le 1$ if and only if dim$(F) \le 1$ and
$[F\colon F ^{q}] \le q$ (see also \cite[Theorem~6.1.8]{GiSz}).
Using the Arason-Baeza theorem and the former part of (1.1), one
computes dim$_{q}(F)$ and CD$(F)$ in the following situation:
\par
\medskip\noindent
(1.2) If $F$ is a field with char$(F) = q > 0$ and $F = F _{\rm
sep}$, such that $[F\colon F ^{q}] = q ^{m}$, for some $m \in
\mathbb{N}$, and $F$ is an extension of an algebraically closed
field $F _{0}$ of transcendence degree $m$, then CD$(F) = {\rm
dim}_{q}(F) = m$.
\par
\medskip
Here it should be pointed out that there exist fields $F _{m}$, $m
\in \mathbb{N}$, of zero characteristic with cd$(\mathcal{G}_{F})
= m < \infty $ and ddim$(F _{m}) = \infty $, for each index $m$;
as shown by Ax, $F _{1}$ can be chosen to be quasifinite (see
\cite{Ax1}). When $m = 2$, one may take as $F _{2}$ any finite
extension of the field $\mathbb{Q}_{p'}$ of $p'$-adic numbers (cf.
\cite{Alemu}, see also Remark \ref{rema2.3}). These results
attract interest in the problem of describing all pairs $(k, \ell
)$ which are equal to $({\rm cd}(\mathcal{G}_{E _{k, \ell }}),
{\rm ddim}(E _{k, \ell }))$, for some field $E _{k, \ell }$. Note
that one may take as $E _{\infty , \infty }$ any purely
transcendental extension $\Phi $ of infinite transcendence degree
over another field $\Phi _{0}$. Since $\Phi $ has a subfield $\Phi
_{0} ^{\prime }$ such that $\Phi _{0} ^{\prime }/\Phi _{0}$ is a
field extension, $\Phi _{0} ^{\prime }$ is $\Phi _{0}$-isomorphic
to $\Phi $, and $\Phi /\Phi _{0} ^{\prime }$ is purely
transcendental of transcendence degree $1$, the equalities
ddim$(\Phi ) = \infty $ and cd$(\mathcal{G}_{\Phi }) = \infty $
can be deduced, by assuming the opposite, from the
Lang-Nagata-Tsen theorem and Galois cohomology, respectively (see
\cite{Na} and \cite[Ch. II, Proposition~11]{S1}). Moreover, it
follows from \cite[Ch. II, Proposition~11]{S1}, and the $\Phi
_{0}$-isomorphism $\Phi \cong \Phi _{0} ^{\prime }$ that
$\mathcal{G}_{\Phi }$ has infinite cohomological $p$-dimensions
cd$_{p}(\mathcal{G}_{\Phi })$, for all $p \in \mathbb{P}$
different from char$(\Phi )$. Hence, by \cite[Theorem~1.15]{KM},
and the (topological) group isomorphism $\mathcal{G}_{\Phi } \cong
\mathcal{G}_{\Phi '}$, where $\Phi ^{\prime }$ is a perfect
closure of $\Phi $ (cf. \cite[Ch. V, Proposition~6.11, and Ch.
VII, Theorem~1.12]{L}), ddim$(\Phi ^{\prime }) = {\rm
cd}(\mathcal{G}_{\Phi '}) = \infty $. This allows us to restrict
ourselves to the case of $(k, \ell ) \neq (\infty , \infty )$,
i.e. $k \neq \infty $. Also, statements (1.1) (and the equalities
cd$(\mathcal{G}_{E}) = {\rm ddim}(E) = 0$, for every algebraically
closed field $E$) imply we may assume further that $k > 0$.
\par
\medskip
The present research solves the stated problem affirmatively, for
each nonzero pair $(k, \ell )$ admissible by Question 1, i.e.
satisfying $1 \le k \le \ell $. It proves that the field $E _{k,
\ell }$ can be chosen to be perfect of any prescribed
characteristic. In order to facilitate our considerations, we
denote by $\mathbb{N} _{\infty }$ and $\overline {\mathbb{P}}$ the
unions $\mathbb{N} \cup \{\infty \}$ and $\mathbb{P} \cup \{0\}$,
respectively, and by $\mathbb{P} _{q}$ the set $\mathbb{P}
\setminus \{q\}$, for each $q \in \overline {\mathbb{P}}$.
\par
\medskip
\section{\bf The main results}
\par
\medskip
The main purpose of this paper is to prove the following:
\par
\medskip
\begin{theo}
\label{theo2.1} For each $q \in \overline {\mathbb{P}}$, there
exist quasifinite fields $F _{m,q}\colon m \in \mathbb{N} _{\infty
}$, of characteristic $q$, such that {\rm ddim}$(F _{m,q}) = m$, for
each $m$.
\end{theo}
\par
\medskip
Theorem \ref{theo2.1} and our next result complement Theorem~1 of
\cite{Ax1} as follows:
\par
\medskip
\begin{theo}
\label{theo2.2} Let $q \in \overline {\mathbb{P}}$ and $(k, \ell)
\in \mathbb{N} \times \mathbb{N} _{\infty }$ be a nonzero pair
admissible by Question 1. Then there exists a perfect field $E
_{k,\ell ; q}$ with {\rm char}$(E _{k, \ell ; q}) = q$, {\rm
ddim}$(E _{k,\ell ; q}) = \ell $ and {\rm cd}$(\mathcal{G}_{E
_{k,\ell ; q}}) = k$.
\end{theo}
\par
\medskip
It is worth mentioning that CD$(E) = {\rm cd}(\mathcal{G}_{E})$,
for every perfect field $E$; in particular, one may write CD$(E
_{k, \ell ; q})$ instead of cd$(\mathcal{G}_{E _{k, \ell ; q}})$
in the statement of Theorem \ref{theo2.2}. In characteristic zero,
the considered equality holds by definition, and for its proof in
case char$(E) = q > 0$, it is sufficient to show that dim$_{q}(E)
= {\rm cd}_{q}(\mathcal{G}_{E}) \le 1$. The inequality
cd$_{q}(\mathcal{G}_{E}) \le 1$ is a well-known result of Galois
cohomology; the same applies to the fact that
cd$_{q}(\mathcal{G}_{E}) = 0$ if and only if $q$ does not divide
the degree $[E ^{\prime }\colon E]$, for any $E ^{\prime } \in
{\rm Fe}(E)$ (cf. \cite[Ch. I, 3.3, and Ch. II, 2.2]{S1}).
Observing further that, for any finite extension $E ^{\prime }$ of
$E$, the group $H _{q} ^{1}(E ^{\prime }) = H
^{1}(\mathcal{G}_{E'}, \mathbb{Z}/q\mathbb{Z})$ is isomorphic to
the group of continuous homomorphisms of the compact group
$\mathcal{G}_{E'}$ into the multiplicative discrete group of
complex $q$-th roots of unity, and using Galois theory and Sylow's
theorem (see \cite[Ch. I, Sect.~6, and Ch. VI]{L}), one concludes
that dim$_{q}(E) = 0$ if and only if cd$_{q}(\mathcal{G}_{E}) =
0$. Thus the equality cd$_{q}(\mathcal{G}_{E}) = {\rm dim}_{q}(E)$
reduces to a consequence of the assertion that dim$_{q}(E) \le 1$.
The assertion itself is true, since $E ^{\prime }$ is a perfect
field, which implies Br$(E ^{\prime })$ does not contain elements
of order $q$ (see \cite[Lemma~9.1.7]{GiSz}) and so proves that $H
_{q} ^{2}(E ^{\prime }) = \{0\}$. \label{k45} One also sees that
cd$(\mathcal{G}_{E}) = \overline {\rm CD}(E)$, where $\overline
{\rm CD}(F)$ is a modification of CD$(F)$, defined for any field
$F$ with char$(F) = q
> 0$, as follows: take as dim$_{q}(F)$ in Definition~1~(ii) the
minimal index $u$ for which $H _{q} ^{u+1}(L) = \{0\}$, for all
finite extensions $L/F$ (cf. \cite[page~715]{ArArBa}); put
dim$_{q}(F) = \infty $ if such $u$ does not exist. In addition, it
is easily verified that: $\overline {\rm CD}(F) \le {\rm CD}(F)$;
$\overline {\rm CD}(F) \le 1$ if and only if dim$(F) \le 1$.
\par
\medskip
\begin{rema}
\label{rema2.3} Note that cd$(\mathcal{G}_{E}) = 2$ and ddim$(E) =
\infty $ in the following two cases: {\rm (i)} $E$ is a finite
extension of the field $\mathbb{Q} _{p}$ of $p$-adic numbers, for
some $p \in \mathbb{P}$; {\rm (ii)} $E$ is a totally imaginary
number field, i.e. a finite extension of the field $\mathbb{Q}$ of
rational numbers, which does not embed in the field $\mathbb{R}$
of real numbers. The equality cd$(\mathcal{G}_{E}) = 2$ is a
well-known result of Galois cohomology (cf. \cite[Ch. II, Sects.
4.3, 4.4]{S1}), and the equality ddim$(E) = \infty $ has been
proved in \cite{Alemu}, for any finite extension $E/\mathbb{Q}
_{p}$ (see also \cite{AK}, for a proof of the fact that
ddim$(\mathbb{Q} _{p}) = \infty $, $p \in \mathbb{P}$). The
validity of the equality ddim$(E) = \infty $ in case {\rm (ii)}
can be deduced from its validity in case {\rm (i)}, by assuming
the opposite and applying \cite[page~379, lemma]{L1} (to suitably
chosen forms without nontrivial zeroes over the completion $E
_{v}$ of $E$ with respect to some discrete valuation $v$).
Algebraic extensions $E _{0}/\mathbb{Q}$ and $E _{p}/\mathbb{Q}
_{p}$, $p \in \mathbb{P}$, such that dim$(E _{p'}) \le 1 < {\rm
ddim}(E _{p'})$, for each $p' \in \overline {\mathbb{P}}$, can be
found in \cite{Ch2}. At the same time, the question of whether
there exist algebraic extensions $E _{0,\ell }/\mathbb{Q}$ and $E
_{p,\ell }/\mathbb{Q} _{p}$, $p \in \mathbb{P}$, with ddim$(E
_{p',\ell }) = \ell $, for each $p' \in \overline {\mathbb{P}}$,
seems to be open, for any integer $\ell \ge 2$.
\end{rema}
\par
\medskip
Theorems \ref{theo2.1} and \ref{theo2.2} are proved in Section 4
and 5, respectively. Our proofs rely on general properties of
fields with Henselian valuations, particularly, of algebraic
extensions of complete discrete valued fields. Preliminaries on
Henselian (valued) fields and other related information used in
the sequel are presented in Section 3. To prove Theorem
\ref{theo2.1}, we show that every algebraically closed field
$\mathbb{F}$ possesses extensions $F _{\infty }$ and $F _{m}$, $m
\in \mathbb{N}$, which are quasifinite fields such that ddim$(F
_{\infty }) = \infty $, and for each index $m$, ddim$(F _{m}) = m$
and $F _{m}$ is a subfield of $F _{\infty }$. The existence of $F
_{\infty }$ has been established constructively by Ax in the case
where char$(\mathbb{F}) = 0$. When char$(\mathbb{F}) = q
> 0$, $F _{\infty }$ is defined by modifying Ax's construction. For
this purpose, we use implicitly the Mel'nikov-Tavgen' theorem
\cite{MT} (via Lemma \ref{lemm3.4}), and the bridge between Galois
theory and the study of tensor products of field extensions,
provided by Cohn's theorem (stated below as Lemma \ref{lemm3.3}).
For the proof of Theorem \ref{theo2.2}, we also need at crucial
points Greenberg's theorem \cite{Gr} (as well as Lemma
\ref{lemm3.2}). The arithmetic ingredient of our proofs, borrowed
from \cite{Ax1}, is based on a version of Vinogradov's theorem on
the ternary Goldbach problem, for certain prime numbers defined as
follows:
\par
\medskip
{\bf Definition 2.} {\it Let $q \in \overline {\mathbb{P}}$, $m$
be an integer $\ge 0$, and $\alpha $ be a real number such that $0
< \alpha < 1$. A number $p \in \mathbb{P}$ is said to be $(m,
\alpha ; q)$-representable if $m = 0$ or $m > 0$ and $p = p _{1} +
p _{2} + p _{3}$, for some $p _{1}, p _{2}, p _{3} \in \mathbb{P}$
with $q < p ^{\alpha } < p _{1} < p _{2} < p _{3}$, which are $(m
- 1, \alpha ; q)$-representable.}
\par
\medskip
For a proof of the following lemma, we refer the reader to
\cite[Lemma~2]{Ax1}.
\par
\medskip
\begin{lemm}
\label{lemm2.4} For each triple $q \in \overline {\mathbb{P}}$, $m
\in \mathbb{N}$, $\alpha \in \mathbb{R}$ with $0 < \alpha < 1$,
there is $c = c(m, \alpha ; q) \in \mathbb{N}$, such that every $p
\in \mathbb{P}$, $p > c$, is $(m, \alpha ; q)$-representable.
\end{lemm}
\par
\medskip
The basic notation, terminology and conventions kept in this paper
are standard and virtually the same as in \cite{TW}, \cite{L} and
\cite{S1}. Throughout, $\mathbb{Z}$ is the additive group of
integers, value groups are written additively, Galois groups are
viewed as profinite with respect to the Krull topology, and by a
profinite group homomorphism, we mean a continuous one.  For any
field $E$, $E ^{\ast }$ is its multiplicative group, $E ^{\ast n}
= \{a ^{n}\colon a \in E ^{\ast }\}$, for each $n \in \mathbb{N}$,
and $\mu (E)$ is the multiplicative group of all roots of unity
lying in $E$. As usual, for any $p \in \mathbb{P}$, $E(p)$ denotes
the maximal $p$-extension of $E$ (in $E _{\rm sep}$), that is, the
compositum of those finite Galois extensions of $E$ in $E _{\rm
sep}$, whose Galois groups are $p$-groups. As in \cite{Ch1},
abrd$_{p}(E)$ stands for the absolute Brauer $p$-dimension of $E$,
defined as the supremum sup$\{{\rm Brd}_{p}(R)\colon R \in {\rm
Fe}(E)\}$. Given a field extension $E ^{\prime }/E$, we write I$(E
^{\prime }/E)$ for the set of intermediate fields of $E ^{\prime
}/E$. When $E ^{\prime }/E$ is a Galois extension, its Galois group
is denoted by $\mathcal{G}(E ^{\prime }/E)$; we say that $E ^{\prime
}/E$ is cyclic if $\mathcal{G}(E ^{\prime }/E)$ is a cyclic group. By
a $\mathbb{Z} _{p}$-extension, we mean a Galois extension $\Psi
^{\prime }/\Psi $ with $\mathcal{G}(\Psi ^{\prime }/\Psi )$
isomorphic to $\mathbb{Z} _{p}$. The value group of any discrete
valued field is assumed to be an ordered subgroup of the additive
group of the field $\mathbb{Q}$; this is done without loss of
generality, in view of \cite[Theorem~15.3.5]{E3}, and the fact that
$\mathbb{Q}$ is a divisible hull of any of its infinite subgroups
(see page \pageref{k99}).
\par
\medskip
\section{\bf Preliminaries on Henselian valuations}
\par
\medskip
For any field $K$ with a (nontrivial) Krull valuation $v$, $O
_{v}(K) = \{a \in K\colon \ v(a) \ge 0\}$ denotes the valuation
ring of $(K, v)$, $M _{v}(K) = \{\mu \in K\colon \ v(\mu ) > 0\}$
the maximal ideal of $O _{v}(K)$, $O _{v}(K) ^{\ast } = \{u \in
K\colon \ v(u) = 0\}$ the multiplicative group of $O _{v}(K)$,
$v(K)$ the value group and $\widehat K = O _{v}(K)/M _{v}(K)$ the
residue field of $(K, v)$, respectively; $\overline {v(K)}$ is a
divisible hull of $v(K)$. The valuation $v$ is said to be
Henselian if it extends uniquely, up-to equivalence, to a
valuation $v _{L}$ on each algebraic extension $L$ of $K$. When
this holds, $(K, v)$ is called a Henselian field. The condition
that $v$ is Henselian has the following two equivalent forms (cf.
\cite[Sect. 18.1]{E3}, or \cite[Theorem~A.12]{TW}):
\par
\medskip\noindent
(3.1) (a) Given a polynomial $f(X) \in O _{v}(K) [X]$ and an element
$a \in O _{v}(K)$, such that $2v(f ^{\prime }(a)) < v(f(a))$, where
$f ^{\prime }$ is the formal derivative of $f$, there is a zero $c
\in O _{v}(K)$ of $f$ satisfying the equality $v(c - a) = v(f(a)/f
^{\prime }(a))$;
\par
(b) For each normal extension $\Omega /K$, $v ^{\prime }(\tau (\mu ))
= v ^{\prime }(\mu )$ whenever  $\mu \in \Omega $, $v ^{\prime }$ is
a valuation of $\Omega $ extending $v$, and $\tau $ is a
$K$-automorphism of $\Omega $.
\par
\medskip
Next we recall some facts concerning the case where $(K, v)$ is a
real-valued field, i.e. $v(K)$ is embeddable as an ordered subgroup
in the additive group $\mathbb{R}$ of real numbers. Fix a completion
$K _{v}$ of $K$ with respect to the topology of $v$, and denote by
$\bar v$ the valuation of $K _{v}$ continuously extending $v$. Then:
\par
\medskip\noindent
(3.2) (a) $(K, v)$ is Henselian if and only if $K$ has no proper
separable (algebraic) extension in $K _{v}$ (cf.
\cite[Corollary~18.3.3]{E3}); in particular, $(K _{v}, \bar v)$ is
Henselian.
\par
(b) The topology of $K _{v}$ as a completion of $K$ is the same as
the one induced by $\bar v$; also, $\bar v(K _{v}) = v(K)$ and
$\widehat K$ equals the residue field of $(K _{v}, \bar v)$ (cf.
\cite[Theorems~9.3.2 and 18.3.1]{E3}).
\par
\medskip
When $v$ is Henselian, so is $v _{L}$, for any algebraic field
extension $L/K$. In this case, we denote by $\widehat L$ the
residue field of $(L, v _{L})$, and put
\par\noindent
$O _{v}(L) = O _{v _{L}}(L)$, $M _{v}(L) = M _{v_{L}}(L)$, $v(L) =
v _{L}(L)$; also, we write $v$ instead of $v _{L}$ when there is
no danger of ambiguity. Clearly, $\widehat L/\widehat K$ is an
algebraic extension and $v(K)$ is an ordered subgroup of $v(L)$,
such that $v(L)/v(K)$ is a torsion group; hence, one may assume
without loss of generality that \label{k99} $v(L)$ is an ordered
subgroup of $\overline {v(K)}$. By Ostrowski's theorem (cf.
\cite[Theorem~17.2.1]{E3}), if $[L\colon K]$ is finite, then it is
divisible by $[\widehat L\colon \widehat K]e(L/K)$, and in case
$[L\colon K] \neq [\widehat L\colon \widehat K]e(L/K)$, the
integer $[L\colon K][\widehat L\colon \widehat K] ^{-1}e(L/K)
^{-1}$ is a power of char$(\widehat K)$ (so char$(\widehat K) \mid
[L\colon K]$); here $e(L/K)$ is the ramification index of $L/K$,
i.e. the index of $v(K)$ in $v(L)$. The extension $L/K$ is called
defectless if $[L\colon K] = [\widehat L\colon \widehat K]e(L/K)$.
When the valuation $v$ is discrete, i.e. $v(K)$ is infinite
cyclic, $L/K$ is defectless in the following two situations:
\par
\medskip\noindent
(3.3) (a) $(K, v)$ is Henselian and $L/K$ is separable (see
\cite[Sect. 17.4]{E3}).
\par
(b) $(K, v)$ is a complete valued field, i.e. $(K, v) = (K _{v},
\bar v)$ (cf. \cite[Ch. XII, Proposition~6.1]{L}). This holds if
$K$ is the Laurent (formal power) series field $K _{0}((X))$ in a
variable $X$ over a field $K _{0}$, and $v$ is the standard
discrete valuation of $K$ inducing on $K _{0}$ the trivial
valuation; then $v(K) = \mathbb{Z}$ and $K _{0}$ is the residue
field of $(K, v)$ (see \cite[Examples~4.2.2 and 9.2.2]{E3}).
\par
\medskip
Assume now that $(K, v)$ is a Henselian field and let $R$ be a
finite extension of $K$. We say that $R/K$ is totally ramified if
$e(R/K) = [R\colon K]$; $R/K$ is called tamely ramified if it is
defectless, $e(R/K)$ is not divisible by char$(\widehat K)$, and
$\widehat R$ is separable over $\widehat K$. The extension $R/K$
is said to be inertial if $[R\colon K] = [\widehat R\colon
\widehat K]$ and $\widehat R/\widehat K$ is separable. Inertial
extensions of $K$ are clearly separable. They have a number of
useful properties, some of which are presented by the following
lemma (for its proof, see \cite[Theorem~A.23]{TW}):
\par
\medskip
\begin{lemm}
\label{lemm3.1}
Let $(K, v)$ be a Henselian field and $K _{\rm ur}$ the compositum of
inertial extensions of $K$ in $K _{\rm sep}$. Then:
\par
{\rm (a)} An inertial extension $R ^{\prime }/K$ is Galois if and
only if so is $\widehat R ^{\prime }/\widehat K$. When this holds,
$\mathcal{G}(R ^{\prime }/K)$ and $\mathcal{G}(\widehat R ^{\prime
}/\widehat K)$ are canonically isomorphic.
\par
{\rm (b)} $v(K _{\rm ur}) = v(K)$, $K _{\rm ur}/K$ is a Galois
extension and $\mathcal{G}(K _{\rm ur}/K) \cong
\mathcal{G}_{\widehat K}$.
\par
{\rm (c)} Finite extensions of $K$ in $K _{\rm ur}$ are inertial, and
the natural mapping of $I(K _{\rm ur}/K)$ into $I(\widehat K _{\rm
sep}/\widehat K)$, by the rule $L \to \widehat L$, is bijective.
\par
{\rm (d)} For each $K _{1} \in {\rm Fe}(K)$, the intersection $K
_{0} = K _{1} \cap K _{\rm ur}$ equals the maximal inertial
extension of $K$ in $K _{1}$; in addition, $\widehat K _{0} =
\widehat K _{1}$.
\end{lemm}
\par
\medskip
Greenberg's theorem \cite{Gr} and the next lemma ensure that if $K =
K _{0}((X))$ and $K _{0}$ is a field with ddim$(K _{0}) < \infty $,
then ddim$(K) = {\rm ddim}(K _{0}) + 1$. They provide some of the
basic tools needed to prove Theorems \ref{theo2.1} and \ref{theo2.2}.
\par
\medskip
\begin{lemm}
\label{lemm3.2} Let $K = K _{0}((X))$ be the Laurent series field
in a variable $X$ over a field $K _{0}$, and let $f _{0}(Y _{1},
\dots , Y _{t})$ be a $K _{0}$-form of degree $d$ in $t$
variables. Fix an algebraic closure $\overline K$ of $K$, denote
by $v$ the standard discrete valuation of $K$ trivial on $K _{0}$,
suppose that $f _{0}$ is without a nontrivial zero over $K _{0}$,
and let $f$ be the $K$-form $\sum _{i=0} ^{d-1} f _{0}(X _{i,1},
\dots , X _{i,t})X ^{i}$, where $X _{i,j}\colon i = 0, \dots , d -
1$; $j = 1, \dots , t$, are algebraically independent variables
over $K$. Then:
\par
{\rm (a)} $f$ is a $K$-form of degree $d$ without a nontrivial
zero over $K$; in addition, $f$ does not possess such a zero over
the perfect closure $\widetilde K$ of $K$ in $\overline K$,
provided that $K _{0}$ is perfect and $d$ is not divisible by {\rm
char}$(K _{0})$;
\par
{\rm (b)} $f _{0}$ does not possess a nontrivial zero over any
totally ramified extension $L$ of $K$ (with respect to $v$); hence,
{\rm ddim}$(K _{0}) \le {\rm ddim}(L)$.
\end{lemm}
\par
\smallskip
\begin{proof}
Let $L$ be a totally ramified extension of $K$. We prove Lemma
\ref{lemm3.2} (b), by assuming that $f _{0}$ has a nontrivial zero
over $L$. Then $f _{0}(Y _{1}, \dots , Y _{t})$ must have a zero
$\alpha = (\alpha _{1}, \dots , \alpha _{t})$, such that $\alpha
_{i} \in O _{v}(L)$, for $i = 1, \dots , t$, and $v _{L}( \alpha
_{i'}) = 0$, for some index $i'$. Since $\widehat L = K _{0}$ and
$f _{0}$ equals its reduction modulo $M _{v}(L)$, this means that
the $t$-tuple $\hat \alpha = (\hat \alpha _{1}, \dots , \hat
\alpha _{t})$ is a nontrivial zero of $f _{0}$ over $L$. The
obtained contradiction proves Lemma \ref{lemm3.2} (b).
\par
The proof of Lemma \ref{lemm3.2} (a) relies on the fact that the
quotient group $v(K)/dv(K)$ is cyclic of order $d$, whose elements
are the cosets $v(X \sp i) + dv(K)$, $i = 0, \dots , d - 1$. Also,
Lemma \ref{lemm3.2} (b) implies that $v(f _{0}(\lambda _{1}, \dots
, \lambda _{t})) \in dv(K)$ whenever $\lambda _{1}, \dots ,
\lambda _{t} \in K$ and $\lambda _{j'} \in K ^{\ast }$, for at
least one index $j'$. This ensures that $\sum _{i=0} ^{d-1} f
_{0}(\lambda _{i,1}, \dots , \lambda _{i,t}) = 0$, where $\lambda
_{i,1}, \dots , \lambda _{i,t} \in K$, for each index $i$, if and
only if all $\lambda _{i,j}$ are equal to zero. Thus the former
part of Lemma \ref{lemm3.2} (a) is proved. For the proof of the
latter one, we may assume that $K _{0}$ is perfect, char$(K _{0})
= q > 0$ and $q \nmid d$. Then it follows from (3.2) and (3.3) (b)
that finite extensions of $K$ in $\widetilde K$ are totally
ramified of $q$-primary degrees. Therefore, the natural embedding
of $K$ into $\widetilde K$ induces canonically a group isomorphism
$v(K)/dv(K) \cong v(\widetilde K)/dv(\widetilde K)$; hence,
$v(\widetilde K)/dv(\widetilde K)$ equals the set $\{v(X ^{i}) +
dv(\widetilde K)\colon i = 0, \dots , d - 1\}$. Since, by Lemma
\ref{lemm3.2} (b), $f _{0}$ does not possess a nontrivial zero
over $\widetilde K$, this allows to prove the latter part of Lemma
\ref{lemm3.2} (a) in the same way as the former one.
\end{proof}
\par
\smallskip
The following lemma (for a proof, see \cite[Theorem~2.2]{Co})
gives a sufficient condition that the tensor product $E _{1}
\otimes _{E _{0}} E _{2}$ is a field, where $E _{i}/E _{0}$, $i =
1, 2$, are field extensions. This lemma considerably simplifies
the presentation of the proof of Theorem \ref{theo2.1}, and its
valuation-theoretic preparation.
\par
\medskip
\begin{lemm}
\label{lemm3.3} Assume that $E _{0}$ is a field and $E _{1}$, $E
_{2}$ are extensions of $E _{0}$ at least one of which is Galois.
Then the tensor product $E _{1} \otimes _{E _{0}} E _{2}$ is a
field if and only if no proper extension of $E _{0}$ is embeddable
as an $E _{0}$-subalgebra in $E _{i}$, $i = 1, 2$. Moreover, if $E
_{0}$ has an extension $F$ such that $E _{i} \in I(F/E _{0})$, $i
= 1, 2$, then $E _{1} \otimes _{E _{0}} E _{2}$ is a field if and
only if $E _{1} \cap E _{2} = E _{0}$; when this holds, $E _{1}
\otimes _{E _{0}} E _{2}$ and the compositum $E _{1}E _{2}$ are
isomorphic as $E _{0}$-algebras.
\end{lemm}
\par
\medskip
Assuming again that $(K, v)$ is a Henselian field, denote by $K
_{\rm tr}$ the compositum of tamely ramified extensions of $K$ in
$K _{\rm sep}$. It is well-known (see \cite[Theorem~A.24]{TW})
that $K _{\rm tr}/K$ is a Galois extension, and finite extensions
of $K$ in $K _{\rm tr}$ are tamely ramified. Our next lemma
presents additional information concerning $K _{\rm tr}/K$. The
stated results are known but we prove them here for convenience of
the reader. The description of some of them relies on the fact
that the quotient $\overline{v(K)}/v(K)$ is an abelian torsion
group, whence it is isomorphic to the direct sum of its
$p$-components $(\overline{v(K)}/v(K)) _{p}$, $p \in \mathbb{P}$.
\par
\medskip
\begin{lemm}
\label{lemm3.4} Let $(K, v)$ be a Henselian field with {\rm
char}$(\widehat K) = q$ and \par\noindent $\mu (\widehat K) = \mu
(\widehat K _{\rm sep})$. Then there is $\Theta \in I(K _{\rm
tr}/K)$ with the following properties:
\par
{\rm (a)} $K _{\rm ur} \cap \Theta = K$, $K _{\rm ur}.\Theta = K
_{\rm tr}$ and $\Theta /K$ is a Galois extension with
$\mathcal{G}(\Theta /K)$ isomorphic to the topological group
product $\prod _{p \in \mathbb{P} _{q}} \Gamma _{p}$, $\Gamma
_{p}$ being the (continuous) character group of the discrete group
$(v(K _{\rm tr})/v(K)) _{p}$, for each index $p$; in particular,
$\mathcal{G}(\Theta /K)$ is abelian.
\par
{\rm (b)} Finite extensions of $K$ in $\Theta $ are totally
ramified, and $\Theta $ equals the compositum of the fields $T
_{p} = \Theta \cap K(p)$, $p \in \mathbb{P} _{q}$; specifically,
if $v$ is discrete, then $T _{p}/K$, $p \in \mathbb{P} _{q}$, are
$\mathbb{Z} _{p}$-extensions.
\par
{\rm (c)} $K _{\rm tr}/K$ is a Galois extension with
$\mathcal{G}(K _{\rm tr}/K) \cong \mathcal{G}(K _{\rm ur}/K) \times
\mathcal{G}(\Theta /K)$.
\par\smallskip
\noindent Moreover, if $q > 0$, then there exists $W \in I(K _{\rm
sep}/K)$, such that $W \cap K _{\rm tr} = K$ and $W.K _{\rm tr} =
K _{\rm sep} = W _{\rm tr}$; also, finite extensions of $K$ in $W$
are of $q$-primary degrees, $W\Theta \cap K _{\rm ur} = K$,
$W\Theta .K _{\rm ur} = (W\Theta ) _{\rm ur} = K _{\rm sep}$ and
$\widehat {W\Theta }/\widehat K$ is a purely inseparable field
extension.
\end{lemm}
\par
\smallskip
\begin{proof}
The existence of a field $\Theta \in I(K _{\rm tr}/K)$ with $K
_{\rm ur} \cap \Theta = K$ and
\par\noindent
$K _{\rm ur}\Theta = K _{\rm tr}$, and in case $q > 0$, the
existence of $W \in I(K _{\rm sep}/K)$, such that $W \cap K _{\rm
tr} = K$ and $W.K _{\rm tr} = K _{\rm sep}$, follow from the
Mel'nikov-Tavgen' theorem \cite{MT} and Galos theory. This,
combined with Lemma \ref{lemm3.1} (d), proves that $K _{\rm tr} =
\Theta _{\rm ur}$, $v(\Theta ) = v(K _{\rm tr})$, and finite
extensions of $K$ in $\Theta $ are totally ramified. Note further
that these extensions are normal with abelian Galois groups.
Indeed, $\mu (\widehat K) = \mu (\widehat K _{\rm sep})$, which
means that $\widehat K$ contains a primitive $\nu $-th root of
unity, for each $\nu \in \mathbb{N}$ not divisible by $q$. This
enables one to obtain the claimed property of the considered
extensions as a consequence of \cite[Proposition~A.22]{TW}.
Moreover, it follows that $\Theta /K$ is Galois,
$\mathcal{G}(\Theta /K)$ is abelian, and because of the equality
$K _{\rm ur} \cap \Delta = K$, there are isomorphisms
$\mathcal{G}(\Theta /K) \cong \mathcal{G}(K _{\rm tr}/K _{\rm
ur})$ and $\mathcal{G}(K _{\rm tr}/K) \cong \mathcal{G}(K _{\rm
ur}/K) \times \mathcal{G}(\Theta /K)$ as profinite groups. Hence,
by Galois theory, the decomposability of $v(K _{\rm tr})/v(K)$
into the direct sum $\oplus _{p \in \mathbb{P} _{q}} (v(K _{\rm
tr}/v(K)) _{p}$, and \cite[Theorem~A.24~(v)]{TW},
$\mathcal{G}(\Theta /K)$ is isomorphic to the product $\prod _{p
\in \mathbb{P} _{q}} \Gamma _{p}$ defined in Lemma \ref{lemm3.4}
(a). Observe that $v(K _{\rm tr})/v(K)$ is divisible. Clearly,
$v(K _{\rm tr})/v(K)$ is an abelian torsion group without an
element of order $q$, so it suffices to prove that $v(K _{\rm tr})
= pv(K _{\rm tr})$, for each $p \in \mathbb{P} _{q}$. Assuming the
opposite and taking an element $\pi \in K _{\rm tr} ^{\ast }$ of
value out of $pv(K _{\rm tr})$, for some $p \in \mathbb{P} _{q}$,
one obtains that the extension $K _{\rm tr} ^{\prime }$ of $K
_{\rm tr}$ generated by some $p$-th root $\pi _{p} \in K _{\rm
sep}$ of $\pi $ must be totally ramified of degree $p$. This
requires that $v(K _{\rm tr} ^{\prime })/v(K)$ is an abelian
torsion group without an element of order divisible by $q$, which
leads to the contradiction that $K _{\rm tr} ^{\prime } = K _{\rm
tr}$. Thus it turns out that $v(K _{\rm tr}) = pv(K _{\rm tr})$,
$p \in \mathbb{P}$, so $v(K _{\rm tr})/v(K)$ is divisible, as
claimed. In addition, it is easily verified that the $p$-group
$(v(K _{\rm tr})/v(K)) _{p}$ is nontrivial, for some $p \in
\mathbb{P} _{q}$, if and only if $v(K) \neq pv(K)$. Suppose now
that $v(K) = \mathbb{Z}$. Then $v(K _{\rm sep}) \cong \mathbb{Q}$,
so it follows from the equality $v(\Theta ) = v(K _{\rm tr})$ and
the cyclicity of finitely-generated subgroups of $\mathbb{Q}$ that
$(v(K _{\rm tr})/v(K)) _{p}$ is a quasi-cyclic $p$-group, for
every $p \in \mathbb{P} _{q}$. This implies the character group of
$(v(K _{\rm tr})/v(K)) _{p}$ is isomorphic to $\mathbb{Z} _{p}$,
which allows to deduce Lemma \ref{lemm3.4} (b) from Galois theory
and Lemma \ref{lemm3.4} (a).
\par
It remains for us to complete the proof of the concluding
assertion of Lemma \ref{lemm3.4}, so we assume that $q > 0$. It
follows from Lemma \ref{lemm3.3} and the equalities $W \cap K
_{\rm tr} = K$, $W.K _{\rm tr} = K _{\rm sep} = W _{\rm tr}$ and
$\Theta \cap K _{\rm ur} = K$, $\Theta .K _{\rm ur} = K _{\rm tr}$
that there exist isomorphisms $K _{\rm sep} \cong W \otimes _{K} K
_{\rm tr}$, $K _{\rm tr} \cong \Theta \otimes _{K} K _{\rm ur}$,
and $W\Theta \cong W \otimes _{K} \Theta $ as $K$-algebras.
Observing also that $K _{\rm sep}$ is $K$-isomorphic
\par\vskip0.064truecm\noindent
to $W \otimes _{K} (\Theta \otimes _{K} K _{\rm ur})$ and $(W
\otimes _{K} \Theta ) \otimes _{K} K _{\rm ur}$, one obtains that
\par\vskip0.064truecm\noindent
$K _{\rm sep} \cong W\Theta \otimes _{K} K _{\rm ur}$, i.e.
$W\Theta \cap K _{\rm ur} = K$ and $W\Theta .K _{\rm ur} = K _{\rm
sep} = (W\Theta ) _{\rm ur}$. Hence, by Lemma \ref{lemm3.1} (d),
$\widehat {W\Theta }/\widehat K$ is a purely inseparable
extension.
\par\vskip0.064truecm
Let now $W _{0}$ be a finite extension of $K$ in $W$. Since $W.K
_{\rm tr} \cong W \otimes _{K} K _{\rm tr}$, it follows that $[W
_{0}\colon K] = [W _{0}.K _{\rm tr}\colon K _{\rm tr}]$ (apply
\cite[Sect. 2, Proposition~c]{P}). This, combined with the fact
that $K _{\rm sep} = K _{\rm tr}(q)$ (cf. \cite{MT}), implies $[W
_{0}\colon K] = [W _{0}K _{\rm tr}\colon K _{\rm tr}]$ is a
$q$-primary number, so Lemma \ref{lemm3.4} is proved.
\end{proof}
\par
\medskip
\section{\bf Proof of Theorem \ref{theo2.1}}
\par
\medskip
The proof of the first main result of this paper relies on the
following sufficient condition for a field $F$ to be quasifinite
with ddim$(F) = \infty $:
\par
\medskip
\begin{prop}
\label{prop4.1} Let $\mathbb{F}$ be an algebraically closed field,
$q = {\rm char}(\mathbb{F})$, and $\mathbb{P} _{q}$ the ordered
set $\{p _{n} \in \mathbb{P}, n \in \mathbb{N}\colon p _{n} \neq
q, p _{n} < p _{n+1},$ for each $n\}$. Assume that $F _{0}$ and $F
_{m}\colon m \in \mathbb{N}$, are perfect fields satisfying the
following conditions:
\par
{\rm (a)} $F _{0} = \mathbb{F}$, provided that $q = 0$; if $q > 0$,
then $\mathcal{G}_{F _{0}} \cong \mathbb{Z} _{q}$ and $F
_{0}/\mathbb{F}$ is a field extension of transcendence degree $1$;
\par
{\rm (b)} For each $m \in \mathbb{N}$, $F _{m}$ is an algebraic
extension of the Laurent series field $F _{m-1}((X _{m}))$,
such that $F _{m-1,{\rm sep}} \otimes _{F _{m-1}} F _{m}$ is a
field and $\mathcal{G}_{F_{m}}$ is isomorphic to the topological
group product $\mathcal{G}_{F_{m-1}} \times \mathbb{Z} _{p_{m}}$.
\par
Then the union $F _{\infty }$ of $F _{n}$, $n \in \mathbb{N}$, is a
quasifinite field with {\rm ddim}$(F) = \infty $.
\end{prop}
\par
\smallskip
\begin{proof}
The fields $F _{m}$, $m \ge 0$, form an ordered set with respect
to inclusion which implies $F _{\infty }$ is a field including $F
_{m}$ as a subfield, for each $m$. Note that $F _{\infty }$ is
perfect. The assertion is obvious if $q = 0$, and in case $q
> 0$, it is easily verified that $F _{\infty } ^{q}$ equals the
union $\cup _{m \in \mathbb{N}} F _{m} ^{q}$. Since the fields $F
_{m}$ are perfect, i.e. $F _{m} ^{q} = F _{m}$, for all integers
$m \ge 0$, this yields $F _{\infty } ^{q} = F _{\infty }$, proving
that $F _{\infty }$ is perfect. We show that $F _{\infty }$ is
quasifinite. Identifying as we can $F _{m,{\rm sep}}$, $m \ge 0$
with their isomorphic copies in $F _{\infty ,{\rm sep}}$, and
taking into account that the polynomial ring $F _{\infty }[X]$ in
an indeterminate $X$ over $F _{\infty }$ equals $\cup _{m \in
\mathbb{N}} F _{m}[X]$, one obtains that $F _{\infty ,{\rm sep}} =
\cup _{m \in \mathbb{N}} F _{m,{\rm sep}}$. Moreover, it follows
from conditions (a) and (b) of Proposition \ref{prop4.1} that $F
_{m-1,{\rm sep}} \cap F _{m} = F _{m-1}$ and $F _{m-1,{\rm sep}}
\otimes _{F _{m-1}} F _{m}$ is $F _{m-1}$-isomorphic to the
compositum $F _{m-1,{\rm sep}}.F _{m}$, for each $m \in
\mathbb{N}$. Arguing by induction on $k$ and observing that $F
_{m-1,{\rm sep}}.F _{m}$ is an $F _{m}$-subalgebra of $F _{m,{\rm
sep}}$, and in case $k \ge 2$, $F _{m-1,{\rm sep}} \otimes _{F
_{m-1}} F _{m-1+k}$ and $(F _{m-1,{\rm sep}} \otimes _{F _{m-1}} F
_{m}) \otimes _{F _{m}} F _{m-1+k}$ are isomorphic $F
_{m-1+k}$-algebras (see \cite[Sect. 9.4, Corollary~a]{P}), one
concludes that $F _{m-1,{\rm sep}} \otimes _{F _{m-1}} F
_{m-1+k}$, $k \in \mathbb{N}$, and $F _{m-1,{\rm sep}} \otimes _{F
_{m-1}} F _{\infty }$ are fields. Hence, by Galois theory (cf.
\cite[Ch. VI, Theorem~1.12]{L}), for each $\mu \in \mathbb{N}$,
every irreducible polynomial $f _{\mu }[X] \in F _{\mu }[X]$ over
$F _{\mu }$ remains irreducible over the fields $F _{\mu '}$, $\mu
< \mu ' \le \infty $; in addition, the Galois groups of $f _{\mu
}(X)$ over $F _{\mu '}$, $\mu \le \mu '$, are isomorphic. It is
now clear from the conditions on $\mathcal{G}_{F _{\mu }}$, $0 \le
\mu < \infty $, that finite extensions of $F _{\infty }$ are
cyclic, and for each $\nu \in \mathbb{N}$, there exists a unique
degree $\nu $ extension of $F _{\infty }$ in $F _{\infty ,{\rm
sep}}$; in other words, the perfect field $F _{\infty }$ is
quasifinite, as claimed.
\par
\smallskip
Our next objective is to prove that ddim$(F _{\infty }) = \infty
$. As a first step towards this goal, we fix an index $\nu \in
\mathbb{N}$ and show that if $g _{\nu }$ is an $F _{\nu }$-form
without a nontrivial zero over $F _{\nu }$, then so is $g _{\nu }$
over $F _{\infty }$. This amounts to proving that $g _{\nu }$ does
not possess a nontrivial zero over $F _{\nu '}$, for any integer
$\nu ' > \nu $. Proceeding by induction on $\nu ' - \nu $, one
reduces our proof to the case where $\nu ' = \nu + 1$. Denote by
$\kappa _{\nu +1}$ the standard valuation of the field $K _{\nu
+1} = F _{\nu }((X _{\nu +1}))$ trivial on $F _{\nu }$, and by $v
_{\nu +1}$ the $\mathbb{Q}$-valued valuation of $F _{\nu +1}$
extending $\kappa _{\nu +1}$. Observing that $F _{\nu ,{\rm
sep}}.K _{\nu +1} = K _{\nu +1,{\rm ur}}$ and $F _{\nu ,{\rm
sep}}.K _{\nu +1} \otimes _{K _{\nu +1}} F _{\nu +1}$ is a field,
one obtains that $F _{\nu ,{\rm sep}}.K _{\nu +1} \cap F _{\nu +1}
= K _{\nu +1}$. This, combined with (3.3) (b) and Lemma
\ref{lemm3.1}, implies finite extensions of $K _{\nu +1}$ in $F
_{\nu +1}$ are totally ramified, so it follows from Lemma
\ref{lemm3.2} (b) that $g _{\nu }$ has no nontrivial zero over $F
_{\nu +1}$. Now our assertion about $g _{\nu }$ is obvious, which
yields ddim$(F _{\infty }) \ge {\rm ddim}(F _{\nu })$, for every $\nu
\in \mathbb{N}$.
\par
\smallskip
It remains to be seen that, for any $m \in \mathbb{N}$, there
exists $\mu (m) \in \mathbb{N}$, such that ddim$(F _{\mu (m)}) >
m$, i.e. $F _{\mu (m)}$ is not of type $C _{m}$. Evidently, there
is a sequence $\beta _{n} \in \mathbb{R}$, $n \in \mathbb{N}$,
satisfying $0 < \beta _{n} < 1$ and $\sum _{j=0} ^{n} \beta _{n}
^{j} > n$, for each index $n$. Therefore, our assertion can be
deduced from the following lemma which applies, by Lemma
\ref{lemm2.4}, to any pair $(\kappa , \beta) \in \mathbb{N} \times
(0, 1)$.
\par
\medskip
\begin{lemm}
\label{lemm4.2} Assume that $F _{\infty }/\mathbb{F}$ is a field
extension, where $\mathbb{F}$ is algebraically closed, {\rm
char}$(\mathbb{F}) = q$, and $F _{\infty }$ is quasifinite; also, let
$F _{0}$ and $F _{n}$, $n \in \mathbb{N}$, be intermediate fields of
$F _{\infty }/\mathbb{F}$ satisfying the conditions of Proposition
\ref{prop4.1}, and such that $\cup _{n \in \mathbb{N}} F _{n} = F
_{\infty }$. Take a pair $(\kappa , \beta ) \in \mathbb{Z} \times
\mathbb{R}$ with $0 < \beta < 1$ and $\kappa \ge 0$, put $\lambda
= \sum _{j=0} ^{\kappa } \beta ^{j}$, and fix some $\mu \in
\mathbb{N}$ so that $p _{\mu } \in \mathbb{P} _{q}$ be $(\beta ,
\kappa ; q)$-representable. Then there exists an $F _{\mu }$-form
$f _{\mu }$ of degree $p _{\mu }$, which depends on at least $p
^{\lambda }$ variables and does not possess a nontrivial zero over
$F _{\mu }$.
\end{lemm}
\par
\medskip
Lemma \ref{lemm4.2} has been proved inductively by Ax in the case
of $q = 0$ (see \cite[Lemma~1]{Ax1}). Here we show that Ax's proof
remains valid for any $q \in \overline{\mathbb{P}}$. Note first
that the field $F _{\nu }$ satisfies condition (b) of Proposition
\ref{prop4.1}, for each $\nu \in \mathbb{N}$, so it follows from
Galois theory that $F _{\nu }$ has an extension $\widetilde F
_{\nu }$ of degree $p _{\nu }$. Hence, the norm form of
$\widetilde F _{\nu }/F _{\nu }$ with respect to any $F _{\nu
}$-basis of $\widetilde F _{\nu }$ is an $F _{\nu }$-form of
degree $p _{\nu }$ in $p _{\nu }$ variables without a nontrivial
zero over $F _{\nu }$, which proves the assertion of Lemma
\ref{lemm4.2} in case $\kappa = 0$. We assume further that $\kappa
\ge 1$, $p _{\mu } = p _{\mu _{1}} + p _{\mu _{2}} + p _{\mu
_{3}}$, where $p _{\mu _{1}}$, $p _{\mu _{2}}$, $p _{\mu _{3}} \in
\mathbb{P} _{q}$ are $(\beta , \kappa - 1; q)$-representable with
$q < p _{\mu } ^{\beta } < p _{\mu _{1}} < p _{\mu _{2}} < p _{\mu
_{3}}$, and the statement of Lemma \ref{lemm4.2} holds for the
$(\beta , \kappa - 1; q)$-representable numbers $p _{\mu _{1}}$,
$p _{\mu _{2}}$ and $p _{\mu _{3}}$. This means that there are $F
_{\mu _{i}}$-forms $f _{\mu _{i}}$ of degree $p _{\mu _{i}}$ ($i =
1, 2, 3$), depending on $d \ge p _{\mu _{1}} ^{\lambda '} > p
_{\mu } ^{\beta \lambda '}$ variables $X _{1}, \dots , X _{d}$,
where $\lambda ' = \sum _{j=0} ^{\kappa - 1} \beta ^{j}$. The
product $\tilde f _{\mu } = f _{\mu _{1}}.f _{\mu _{2}}.f _{\mu
_{3}}$ is clearly an $F _{\mu -1}$-form of degree $p _{\mu }$ in
variables $X _{1}, \dots , X _{d}$. Using Galois theory, Lemma
\ref{lemm3.4}, and condition (b) of Proposition \ref{prop4.1}, one
proves that $v _{\mu }(F _{\mu }) \neq p _{\mu }v _{\mu }(F _{\mu
})$, $v _{\mu }$ being the $\mathbb{Q}$-valued valuation of $F
_{\mu }$ extending the standard valuation of $F _{\mu -1}((X _{\mu
}))$. Fix an element $\pi \in F _{\mu } ^{\ast }$ with $v _{\mu
}(\pi ) \notin p _{\mu }v _{\mu }(F _{\mu })$, and let $f _{\mu }$
be the $F _{\mu }$-form $\sum _{u=0} ^{p _{\mu }-1} \tilde f _{\mu
}(X _{u,1}, \dots , X _{u,d})\pi ^{u}$ (in $dp _{\mu } > p _{\mu }
^{1+\beta \lambda '} = p _{\mu } ^{\lambda }$ variables). It
follows from the above-noted properties of $F _{\nu }$-forms, for
$\nu \in \mathbb{N}$, that $\tilde f _{\mu }$ does not possess a
nontrivial zero over $F _{\mu -1}$. Since the fulfillment of
condition (b) guarantees that finite extensions of $F _{\mu -1}((X
_{\mu }))$ in $F _{\mu }$ are totally ramified, this allows to
deduce from Lemma \ref{lemm3.2}~(b) that $\tilde f _{\mu }$ is
without a nontrivial zero over $F _{\mu }$. Therefore, one obtains
as in the proof of Lemma \ref{lemm3.2}~(a) that $f _{\mu }$ has
the properties claimed by Lemma \ref{lemm4.2}, so Proposition
\ref{prop4.1} is proved.
\end{proof}
\par
\smallskip
Next we show that every algebraically closed field $\mathbb{F}$
has extensions subject to the restrictions of Proposition
\ref{prop4.1}. We begin with the following lemma.
\par
\smallskip
\begin{lemm}
\label{lemm4.3}
Let $\mathbb{F}$ be an algebraically closed field of characteristic $q
> 0$, $\Pi $ a nonempty subset of $\mathbb{P}$, and $\mathbb{F} _{1}
= \mathbb{F}(X)$ be the rational function field over $\mathbb{F}$ in a
variable $X$. Then there exists an algebraic extension $F
_{0}/\mathbb{F} _{1}$, such that $F _{0}$ is a perfect field and
$\mathcal{G}_{F _{0}} \cong \prod _{p \in \Pi } \mathbb{Z} _{p}$.
\end{lemm}
\par
\smallskip
\begin{proof}
Let $\overline {\mathbb{F}} _{1}$ be an algebraic closure of
$\mathbb{F} _{1,{\rm sep}}$, and for each $p \in \mathbb{P}$, let
$L _{p}$ be an extension of $\mathbb{F} _{1}$ in $\overline
{\mathbb{F}} _{1}$ obtained by adjunction of a root of the
polynomial $Y ^{p} - \delta _{p}Y - X ^{-1}$, where $\delta _{p} =
1$ if $p = q$, and $\delta _{p} = 0$, otherwise. The field
$\mathbb{F} _{1}$ has a valuation $v$ with $v(\mathbb{F} _{1}) =
\mathbb{Z}$ and $v(X) = 1$ (see \cite[Example~4.1.3]{E3}); this
implies $v$ is trivial on $\mathbb{F}$ and $L _{p}/\mathbb{F}
_{1}$, $p \in \mathbb{P}$, are cyclic field extensions of degree
$p$. As $\mathbb{F} _{1}$ contains a primitive $h$-th root of
unity, for each $h \in \mathbb{N}$ not divisible by $q$, it
follows from Kummer theory and Witt's lemma (cf. \cite[Sect.~15,
Lemma~2]{Dr1}) that $\mathbb{F} _{1}$ has $\mathbb{Z}
_{p}$-extensions $\Gamma _{p}$ in $\overline {\mathbb{F}} _{1}$,
$p \in \mathbb{P}$, such that $L _{p} \in I(\Gamma _{p}/\mathbb{F}
_{1})$, for each $p$. Note also that the compositum $\Gamma $ of
the fields
\par\noindent
$\Gamma _{p}$, $p \in \Pi $, is a Galois extension of $\mathbb{F}
_{1}$ with $\mathcal{G}(\Gamma /\mathbb{F} _{1}) \cong \prod _{p
\in \Pi } \mathbb{Z} _{p}$. Hence, $\mathcal{G}(\Gamma /\mathbb{F}
_{1})$ is a projective profinite group, in the sense of \cite[Ch.
I, 5.9]{S1}, which allows to obtain from Galois theory that there
exists $\widetilde F _{0} \in I(\mathbb{F} _{1,{\rm
sep}}/\mathbb{F} _{1})$ satisfying $\Gamma \cap \widetilde F _{0}
= \mathbb{F} _{1}$ and $\Gamma .\widetilde F _{0} = \mathbb{F}
_{1,{\rm sep}}$. It is proved similarly that if $F _{0}$ is the
perfect closure of $\widetilde F _{0}$ in $\overline {\mathbb{F}}
_{1}$, then $\Gamma \cap F _{0} = \mathbb{F} _{1}$ and $\Gamma F
_{0} = \overline {\mathbb{F}} _{1}$; in particular, $F _{0}$ is a
perfect field with $\mathcal{G}_{F _{0}} \cong \mathcal{G}(\Gamma
/\mathbb{F} _{1}) \cong \prod _{p \in \Pi } \mathbb{Z} _{p}$.
\end{proof}
\par
\smallskip
Let now $\mathbb{F}$ be an algebraically closed field, $F
_{0}/\mathbb{F}$ a field extension satisfying condition (a) of
Proposition \ref{prop4.1}, and let $F _{n}$, $n \in \mathbb{N}$, be
perfect fields with Henselian $\mathbb{Q}$-valued valuations $v
_{n}$, defined inductively as follows:
\par
\medskip\noindent
(4.1) For each $n \in \mathbb{N}$, $F _{n}$ is an extension of the
Laurent series field
\par\noindent
$F _{n-1}((X _{n})) \colon = K _{n}$ in an algebraic closure
$\overline{K} _{n}$ of $F _{n-1,{\rm sep}}.K _{n}$, such that:
\par
(a) $(F _{n}, v _{n})/(K _{n}, \kappa _{n})$ is a valued field
extension, where $\kappa _{n}$ is the standard $\mathbb{Z}$-valued
valuation of $K _{n}$ trivial on $F _{n-1}$, and $v _{n}$ is the
valuation of $F _{n}$ extending $\kappa _{n}$; hence, $(F _{n}, v
_{n})$ is Henselian with $\mathbb{Z} \subseteq v _{n}(F _{n})
\subseteq \mathbb{Q}$;
\par
(b) $F _{n}$ contains as a subfield a separable extension $W _{n}$ of
$K _{n}$ in $\overline{K} _{n}$, such that $W _{n} \cap K _{n,{\rm
tr}} = K _{n}$ and $W _{n}.K _{n,{\rm tr}} = W _{n,{\rm tr}}$ equals
the separable closure $K _{n,{\rm sep}}$ of $K _{n}$ in $\overline{K}
_{n}$; in particular, if $q = 0$, then $W _{n} = K _{n}$;
\par
(c) $F _{n}$ is the extension of $W _{n}$ generated by the $h
_{n}$-th roots of $X _{n}$ in $\overline{K} _{n}$, where $h _{n}$
runs across the set of positive integers not divisible by $p
_{n}$.
\par
\medskip\noindent
We show that the fields $F _{n}$, $n \in \mathbb{N}$, satisfy
condition (b) of Proposition \ref{prop4.1}. Proceeding by
induction on $n$, one obtains that it is sufficient to prove our
assertion, for a fixed index $n$, under the hypothesis that $F
_{n-1}$ is perfect and satisfies condition (a) or (b) of
Proposition \ref{prop4.1} depending on whether or not $n = 1$. Let
$T _{n,p}$ be the extension of $K _{n}$ in $\overline{K} _{n}$
generated by all roots of $X _{n}$ of $p$-primary degrees, for
each $p \in \mathbb{P} _{q}$, and let $T _{n}$ be the compositum
of the fields $T _{n,p}$, $p \in \mathbb{P} _{q}$. Clearly, $X
_{n} \notin K _{n} ^{\ast p}$, for any $p \in \mathbb{P}$, so it
follows from Kummer theory and (3.3) that, in case $p \neq q$, $T
_{n,p}$ is a $\mathbb{Z} _{p}$-extension of $K _{n}$ in $K
_{n,{\rm tr}}$, and finite extensions of $K _{n}$ in $T _{n,p}$
are totally ramified of $p$-primary degrees. As $K _{n,{\rm ur}} =
F _{n-1,{\rm sep}}.K _{n}$, these facts and Lemma \ref{lemm3.4}
(a) show that $F _{n-1,{\rm sep}}.K _{n} \cap T _{n} = K _{n}$ and
$F _{n-1,{\rm sep}}T _{n} = K _{n,{\rm tr}}$. On the other hand,
if $q > 0$, then $K _{n} ^{q} = F _{n-1}((X _{n} ^{q}))$, which
implies $[K _{n}\colon K _{n} ^{q}] = q$ and $T _{n,q}$ is the
perfect closure of $K _{n}$ in $\overline{K} _{n}$. Put
$\widetilde W _{n} = W _{n}T _{n,q}$ if $q > 0$, and $\widetilde W
_{n} = K _{n}$ if $q = 0$. It is easily verified that $\widetilde
W _{n,{\rm ur}} = F _{n-1,{\rm sep}}.\widetilde W _{n}$. Using
(4.1), one obtains that $F _{n} = \widetilde W _{n}.\Phi _{n}$,
$\Phi _{n}$ being the compositum of the fields $T _{n,p}$, for $p
\in \mathbb{P} \setminus \{q, p _{n}\}$. Therefore, $F _{n}$ is
perfect, $F _{n,{\rm tr}} = \overline{K} _{n}$, and in case $q =
0$, we have $F _{n} = \Phi _{n}$, $F _{n-1,{\rm sep}}T _{n, p
_{n}} \cap F _{n} = K _{n}$, $F _{n-1,{\rm sep}}T _{n, p _{n}}.F
_{n} = F _{n,{\rm sep}}$. When $q
> 0$, it follows that $\widetilde W _{n} \cap K _{n,{\rm tr}} = K
_{n}$, $\widetilde W _{n}.K _{n,{\rm tr}} = \overline{K} _{n}$,
and $\widetilde W _{n}\Phi _{n} = F _{n}$. Putting
\par\noindent
$\Theta _{n,p} = \widetilde W _{n}T _{n,p}$, for each $p \in
\mathbb{P} _{q}$, one obtains that $F _{n}$ equals the compositum
of the fields $\Theta _{n,p}$, $p \in \mathbb{P} \setminus \{q, p
_{n}\}$. Since, by (3.3), (4.1) (b) and Lemma \ref{lemm3.4},
finite extensions of $K _{n}$ in $\widetilde W _{n}$ are totally
ramified of $q$-primary degrees, this allows to deduce from Galois
theory and the definition of $T _{n,p}$, $p \in \mathbb{P} _{q}$,
that $\Theta _{n,p}/\widetilde W _{n}$ is a $\mathbb{Z}
_{p}$-extension whose finite subextensions are totally ramified,
for each $p \neq q$. One also sees that finite extensions of
$\widetilde W _{n}$ in $F _{n}$ are tamely and totally ramified of
degrees not divisible by $p _{n}$. Similarly, it is proved that,
for any finite extension $\widetilde W _{n} ^{\prime }$ of
$\widetilde W _{n}$ in $F _{n-1,{\rm sep}}T _{n, p
_{n}}.\widetilde W _{n}$, $e(\widetilde W _{n} ^{\prime
}/\widetilde W _{n})$ is a $p _{n}$-primary number. These
observations show that $F _{n-1,{\rm sep}}\Theta _{n, p _{n}}$ is
a Galois extension of $\widetilde W _{n}$, such that $F _{n-1,{\rm
sep}}\Theta _{n, p _{n}} \cap F _{n} = \widetilde W _{n}$, so it
follows from Galois theory that $\mathcal{G}_{F _{n}} \cong
\mathcal{G}(F _{n-1,{\rm sep}}\Theta _{n, p _{n}}/\widetilde W
_{n}) \cong \mathcal{G}_{F _{n-1}} \times \mathbb{Z} _{p _{n}}$.
In view of our hypothesis on $F _{n-1}$, this completes the proof
of the assertion that the fields $F _{n}$, $n \in \mathbb{N}$,
defined by (4.1) satisfy condition (b) of Proposition
\ref{prop4.1}. Hence, the union $F _{\infty } = \cup _{n \in
\mathbb{N}} F _{n}$ is a quasifinite field with ddim$(F _{\infty
}) = \infty $. As $\mathbb{F}$ is an arbitrary algebraically
closed field, the existence of quasifinite fields $F _{\infty
,q}$, $q \in \overline {\mathbb{P}}$, such that char$(F _{\infty
,q}) = q$ and ddim$(F _{\infty ,q}) = \infty $, for each $q$, is
now obvious.
\par
\smallskip
We turn to the proof of Theorem \ref{theo2.1} in the case where $m
< \infty $. Fix an algebraically closed field $\mathbb{F}$ of
characteristic $q$ as well as an extension $F _{0}$ of
$\mathbb{F}$ satisfying condition (a) of Proposition
\ref{prop4.1}, and a sequence of fields $K _{n}$, $\overline K
_{n}$, $W _{n}$ and $F _{n}$, $n \in \mathbb{N}$, defined in
agreement with (4.1). If $m = 1$ and $q = 0$, then one may put $F
_{m,q} = K _{1}$, and in case $m = 1$ and $q > 0$, $F _{m,q}$ may
be defined by applying Lemma \ref{lemm4.3} to the set $\Pi =
\mathbb{P}$. Henceforth, we assume that $m \ge 2$. It follows from
Greenberg's theorem (and Tsen's theorem, for $n = 0$), that
ddim$(F _{n}) \le n + 1$, for every $n \in \mathbb{N}$. Let
char$(F _{0}) = q$ and $\mu $ be the minimal integer for which
ddim$(F _{\mu }) \ge m$. Using Greenberg's theorem and the
closeness of the class of $C _{n}$-fields under the formation of
algebraic extensions, for each $n \in \mathbb{N}$, one obtains
that ddim$(\Psi _{\mu }) = m$ whenever \par\noindent $\Psi _{\mu }
\in I(F _{\mu }/K _{\mu })$. Note that $\Psi _{\mu }$ can be
chosen to be a quasifinite field. If $\mu = 1$ (which requires
that $q > 0$), then one may take as $\Psi _{\mu }$ the perfect
closure of $W _{\mu }$ in $\overline{F} _{\mu }$. Suppose further
that $\mu \ge 2$. It can be deduced from Galois theory and Lemma
\ref{lemm3.4} that $\mathcal{G}_{W _{\mu }}$ is isomorphic to the
profinite groups $\mathcal{G}(K _{\mu ,{\rm tr}}/K _{\mu })$ and
$\mathcal{G}_{F _{\mu -1}} \times \prod _{p \in \mathbb{P} _{q}}
\mathbb{Z} _{p}$. This implies the existence of a Galois extension
$\Psi _{\mu } ^{\prime }$ of $W _{\mu }$ in $F _{\mu }$ with
$\mathcal{G}(\Psi _{\mu } ^{\prime }/W _{\mu }) \cong \prod _{j=1}
^{\mu -1} \mathbb{Z} _{p _{j}}$. Since $F _{0}$ is algebraically
closed if $q = 0$, and there are isomorphisms $\mathcal{G}_{F
_{0}} \cong \mathbb{Z} _{q}$ and $\mathcal{G}_{F _{\mu -1}} \cong
\mathbb{Z} _{q} \times \prod _{j=1} ^{\mu -1} \mathbb{Z} _{p
_{j}}$ in case $q > 0$, it is easy to see that $\mathcal{G}_{\Psi
_{\mu }'} \cong \prod _{p \in \mathbb{P}} \mathbb{Z} _{p}$. Let
finally $\Psi _{\mu }$ be the perfect closure of $\Psi _{\mu }
^{\prime }$ in $F _{\mu ,{\rm sep}} = \overline K _{\mu }$. Then
it follows from Galois theory and \cite[Ch. V,
Proposition~6.11]{L}, that $\mathcal{G}_{\Psi _{\mu }} \cong
\mathcal{G}_{\Psi _{\mu }'}$, which means that $\Psi _{\mu }$ is a
quasifinite field. Taking into account that $F _{\mu }$ is
perfect, one concludes that $\Psi _{\mu } \in I(F _{\mu }/W _{\mu
})$, so Theorem \ref{theo2.1} is proved.
\par
\smallskip
\begin{coro}
\label{coro4.4}
Let $F _{\infty }$ be a quasifinite field defined in accordance with
Proposition \ref{prop4.1}. Then {\rm ddim}$(F _{\infty } ^{\prime })
= \infty $, for every finite extension $F _{\infty } ^{\prime }/F
_{\infty }$.
\end{coro}
\par
\smallskip
\begin{proof}
By definition, $F _{\infty }$ possesses subfields $\mathbb{F}$,
$F _{0}$ and $F _{n}$, $n \in \mathbb{N}$, such that $\mathbb{F}$ is
algebraically closed, $F _{0}$ and $F _{n}$, $n \in \mathbb{N}$,
satisfy conditions (a) and (b) of Proposition \ref{prop4.1},
respectively, and $F _{\infty } = \cup _{n \in \mathbb{N}} F _{n}$.
We first show that the
fields $K _{n} = F _{n-1}((X _{n}))$, $n \in \mathbb{N}$, given in
condition (b), can be chosen so that $p _{n} \nmid [\Lambda
_{n}\colon K _{n}]$, for any finite extension $\Lambda _{n}$ of $K
_{n}$ in $F _{n}$ (and any $n \in \mathbb{N}$). Proceeding by
induction on $n$, one obtains that there exist elements $Y _{n}
\in F _{n}$, $n \in \mathbb{N}$, such that $Y _{n}$ is a root of
$X _{n}$ of $p _{n}$-primary degree and the degrees of finite
extensions of $K _{n}(Y _{n})$ in $F _{n}$ are not divisible by $p
_{n}$, for any $n$. Since $K _{n}(Y _{n}) = F _{n-1}((Y _{n}))$,
$n \in \mathbb{N}$, this yields the desired reduction.
\par
Take a finite extension $F _{\infty } ^{\prime }$ of $F _{\infty
}$ in $F _{\infty ,{\rm sep}}$, put $\theta = [F _{\infty }
^{\prime }\colon F _{\infty }]$, and for each $n \in \mathbb{N}$,
let $\theta _{n}$ be the greatest $p _{n}$-primary divisor of
$\theta $; put $\theta _{0} = 1$ if $q = 0$, and let $\theta _{0}$
be the greatest $q$-primary divisor of $\theta $ in case $q
> 0$. Fix a primitive element $\lambda $ of $F _{\infty } ^{\prime
}/F _{\infty }$, and denote by $f(X)$ the minimal polynomial of
$\lambda $ over $F _{\infty }$. Clearly, the coefficients of $f(X)$
lie in $F _{m}$, for some $m < \infty $. Consider the fields $F
_{n} ^{\prime }$, $n \ge 0$, defined as follows: $F _{n} ^{\prime
} = F _{n}(\lambda )$ if $n \ge m$; $F _{n} ^{\prime } = F
_{n,{\rm sep}} \cap F _{m} ^{\prime }$ if $n < m$. It follows from
condition (b) that $[F _{n} ^{\prime }\colon F _{n}] \le \theta $,
for all $n \ge 0$, and equality holds if $n \ge m$. Since
$\mathbb{Z} _{p}$ is isomorphic to its open subgroups, for each $p
\in \mathbb{P}$, this means that $F _{0} ^{\prime }$ satisfies
condition (a) of Proposition \ref{prop4.1}. Similarly, it follows
that $\mathcal{G}_{F _{n}'} \cong \mathcal{G}_{F _{n}}$, for every
$n \in \mathbb{N}$.
\par
Note further that, by Lemma \ref{lemm3.3} and condition (b), $F
_{n-1,{\rm sep}} \cap F _{n} = F _{n-1}$ and $F _{n-1,{\rm sep}}
\otimes _{F _{n-1}} F _{n} \cong F _{n-1,{\rm sep}}.F _{n}$ as an $F
_{n-1}$-algebra. These observations and our choice of the fields $K
_{n}$, $n \in \mathbb{N}$, enable one to deduce from Galois
theory that $[F _{n} ^{\prime }\colon F _{n}] = \prod _{j=0} ^{n}
\theta _{j}$, for each $n < m$, and when $n > 0$,
\par\vskip0.04truecm\noindent
$F _{n} ^{\prime } = F _{n-1} ^{\prime }.F _{n}(Z _{n})$ and $F
_{n-1} ^{\prime }.K _{n}(Z _{n}) = F _{n-1} ^{\prime }((Z _{n}))$,
where $Z _{n} \in F _{n} ^{\prime }$ is a
\par\vskip0.04truecm\noindent
$\theta _{n}$-th root of $X _{n}$. This ensures that $F _{n} ^{\prime
} = F _{n}.F _{n-1} ^{\prime }((Z _{n}))$ and finite extensions of $F
_{n-1} ^{\prime }((Z _{n}))$ in $F _{n} ^{\prime }$ are totally
ramified. It is therefore clear from (3.3) (b) and the equality $F
_{n-1,{\rm sep}}.F _{n-1} ^{\prime }((Z _{n})) = F _{n-1} ^{\prime
}((Z _{n})) _{\rm ur}$ that
\par\vskip0.057truecm\noindent
$F _{n-1,{\rm sep}}.F _{n-1} ^{\prime }((Z _{n})) \cap F
_{n} ^{\prime } = F _{n-1} ^{\prime }((Z _{n}))$, which in turn
implies
\par\vskip0.057truecm\noindent
$F _{n-1,{\rm sep}} \cap F _{n} ^{\prime } = F _{n-1,{\rm sep}} \cap
F _{n-1} ^{\prime }((Z _{n})) = F _{n-1} ^{\prime }$. Applying now
Lemma \ref{lemm3.3},
\par\vskip0.053truecm\noindent
one concludes that $F _{n-1,{\rm sep}}.F _{n} ^{\prime } \cong F
_{n-1,{\rm sep}} \otimes _{F _{n-1}'} F _{n} ^{\prime }$ as $F _{n-1}
^{\prime }$-algebras, for every $n > 0$. Thus the extensions $F _{n}
^{\prime }/F _{n-1} ^{\prime }((Z _{n}))$, $n \in \mathbb{N}$ (where
$Z _{n} = X _{n}$ if $n > m$), satisfy condition (b) of Proposition
\ref{prop4.1}, so ddim$(F _{\infty } ^{\prime }) = \infty $, as
required.
\end{proof}
\par
\smallskip
\begin{rema}
\label{rema4.5} Let $F _{m,q}$, $(m, q) \in \mathbb{N} _{\infty }
\times \mathbb{P}$, be quasifinite fields such that char$(F
_{m,q}) = q$ and ddim$(F _{m,q}) = m$, for each pair $(m, q)$.
Then, by Witt's theorem (cf. \cite[Theorem~12.4.1]{E3}), there are
complete discrete valued fields $(K _{m,q}, v _{m,q})$ with
$char(K _{m,q}) = 0$ and $\widehat K _{m,q} = F _{m,q}$, for all
$(m, q) \in \mathbb{N} _{\infty } \times \mathbb{P}$. Also, it
follows from Lemma \ref{lemm3.4} that, for each $(m, q)$, there
exists an extension $T _{m,q}$ of $K _{m,q}$ in $K _{m,q,{\rm
sep}}$, such that $K _{m,q,{\rm ur}}.T _{m,q} = K _{m,q,{\rm
sep}}$ and $K _{m,q,{\rm ur}} \cap T _{m,q} = K _{m,q}$. This
ensures that $T _{m,q,{\rm ur}} = K _{m,q,{\rm sep}}$ and
$\widehat T _{m,q} = F _{m,q}$, so Lemma \ref{lemm3.1} (b) implies
$T _{m,q}$ is a quasifinite field, ddim$(T _{m,q}) \ge m$ (and
equality holds if $m = \infty $). Note that $T _{m,q}$, $(m,q) \in
\mathbb{N} _{\infty } \times \mathbb{P}$, are pairwise
non-isomorphic fields. Since the valuation of $T _{m,q}$ extending
$v _{m,q}$ is Henselian with $v _{m,q}(T _{m,q}) = \mathbb{Q}$,
for every $(m, q)$, this can be proved by assuming the opposite,
and using the non-existence of a field with a pair of Henselian
real-valued valuations whose residue fields are quasifinite and
non-isomorphic. The noted fact follows from the validity of
Schmidt's Uniqueness Theorem in the case of Henselian real-valued
valuations on a field (see \cite[Corollary~21.1.2]{E3}). It would
be of interest to know whether ddim$(T _{m,q}) = m$, for every
$(m, q) \in \mathbb{N} \times \mathbb{P}$.
\end{rema}
\par
\smallskip
\section{\bf Perfect fields of characteristic $q \in \overline
{\mathbb{P}}$ and prescribed Galois cohomological and Diophantine
dimensions}
\par
\medskip
In this Section we present a proof of Theorem \ref{theo2.2}. Our
starting point is the following lemma which proves Theorem
\ref{theo2.2} in case $q = 0$:
\par
\medskip
\begin{lemm}
\label{lemm5.1} Let $\ell \in \mathbb{N} _{\infty }$ and $k$
be an integer satisfying $1 \le k \le \ell $. Then there exists a
field $E _{k,\ell }$ with {\rm ddim}$(E _{k,\ell }) = \ell $ and {\rm
cd}$(E _{k,\ell }) = k$. Furthermore, for each $q \in \overline
{\mathbb{P}}$, $E _{k,\ell }$ can be chosen so that {\rm
char}$(E _{k,\ell }) = q$.
\end{lemm}
\par
\medskip
\begin{proof}
Our assertion is contained in Theorem \ref{theo2.1} in case $k =
1$, so we assume that $k \ge 2$. If $\ell \in \mathbb{N}$, then it
follows from Galois cohomology (see \cite[pages~1219-1220]{Ax1},
and \cite[Ch. II, 2.2 and 4.3]{S1}), Greenberg's theorem and Lemma
\ref{lemm3.2} that one may take as $E _{k,\ell }$ the iterated
Laurent formal power series field $F _{\ell -k+1,q}((X _{1}))
\dots ((X _{k-1}))$, where $F _{\ell -k+1,q}$ has the properties
required by Theorem \ref{theo2.1}, for $m = \ell - k + 1$ and any
$q \in \overline {\mathbb{P}}$. Similarly, if $\ell = \infty $ and
$F _{\infty ,q}$ is a quasifinite field with char$(F _{\infty ,
q}) = q$ and ddim$(F _{\infty , q}) = \infty $, then one may put
$E _{k,\infty } = F _{\infty ,q}((X _{1})) \dots ((X _{k-1}))$.
\end{proof}
\par
\medskip
Theorem \ref{theo2.1} and Lemma \ref{lemm5.1} allow to assume in
the rest of the proof of Theorem \ref{theo2.2} that $q > 0$ and $k
\ge 2$. Retaining notation as in the proof of Lemma \ref{lemm5.1},
denote by $E _{k, \ell } ^{\prime }$ the perfect closure of $E
_{k, \ell }$ in its algebraic closure $\overline E _{k, \ell }$.
It follows from Galois theory and \cite[Ch. V,
Proposition~6.11]{L}, that $\mathcal{G}_{E _{k, \ell }} \cong
\mathcal{G}_{E _{k, \ell }'}$, so it suffices to prove that the
quasifinite constant field of $E _{k, \ell }$ can be chosen so
that ddim$(E _{k, \ell } ^{\prime }) = \ell $. Suppose first that
$\ell = \infty $ and $F _{\infty , q}$ is a quasifinite field with
char$(F _{\infty , q}) = q$ and ddim$(F _{\infty , q}) = \infty $,
defined as in Proposition \ref{prop4.1}. Then there exist $q _{n}
\in \mathbb{P} _{q}$, $n \in \mathbb{N}$, such that for each $n$,
there is an $F _{\infty , q}$-form $t _{n}$ of degree $q _{n}$ in
at least $q _{n} ^{n}$ variables, which does not possess a
nontrivial zero over $F _{\infty , q}$. Proceeding by induction on
$k - 1$ and using Lemma \ref{lemm3.2}, one obtains that $t _{n}$
is without a nontrivial zero over $E _{k, \infty } ^{\prime }$.
This yields ddim$(E _{k, \infty } ^{\prime }) = \infty $, which
completes our proof in case $\ell = \infty $.
\par
\medskip
To prove Theorem \ref{theo2.2} in case $\ell < \infty $ we need
the following two lemmas.
\par
\medskip
\begin{lemm}
\label{lemm5.2} Let $\mathbb{F}$ be an algebraically closed field,
$q = {\rm char}(\mathbb{F})$, $m$ an integer $\ge 2$, and $\Pi
_{i}\colon i = 0, \dots , m - 1$, be nonempty subsets of
$\mathbb{P}$, such that $q \in \Pi _{0}$ in case $q > 0$, $\cup
_{i=0} ^{m-1} \Pi _{i} = \mathbb{P}$, and $\Pi _{i'} \cap \Pi _{i''}
= \emptyset $, provided that $0 \le i' < i'' \le m - 1$. Then there
exist perfect fields $F _{i}\colon i = 0, \dots , m - 1$, with the
following properties:
\par
{\rm (a)} $\mathcal{G}_{F _{0}} \cong \prod _{p \in \Pi _{0}}
\mathbb{Z} _{p}$ and $F _{0}/\mathbb{F}$ is an extension of
transcendence degree $1$;
\par
{\rm (b)} For each $i > 0$, $F _{i}$ is an algebraic extension of the
Laurent series field $F _{i-1}((X _{i}))$, such that $F
_{i-1,{\rm sep}} \otimes _{F _{i-1}} F _{i}$ is a field and
$\mathcal{G}_{F_{i}} \cong \mathcal{G}_{F_{i-1}} \times \prod _{q
_{i} \in \Pi _{i}} \mathbb{Z} _{q_{i}}$; in particular, $F _{m-1}$ is
a quasifinite field with {\rm ddim}$(F _{m-1}) \le m$.
\end{lemm}
\par
\medskip
\begin{proof}
The existence of $F _{0}$ follows from Lemma \ref{lemm4.3}, and the
equality
\par\noindent
ddim$(F _{0}) = 1$ is implied by Tsen's theorem and the fact that $F
_{0} \neq F _{0,{\rm sep}}$. Suppose that $i > 0$ and
$F _{i-1}$ has been defined in accordance with Lemma
\ref{lemm5.2} (a) or (b) depending on whether or not $i = 1$. Put $K
_{i} = F _{i-1}((X _{i}))$, fix an algebraic closure $\overline K
_{i}$ of $K _{i,{\rm sep}}$, and denote by $\kappa _{i}$ the standard
discrete valuation of $K _{i}$ trivial on $F _{i-1}$. Considering the
Henselian field $(K _{i}, \kappa _{i})$, let $W _{i}$ be an extension
of $K _{i}$ in $\overline K _{i}$, such that $W _{i}.K _{i,{\rm tr}}
= \overline K _{i} = W _{i,{\rm sep}}$ and $W _{i} \cap K _{i,{\rm
tr}} = K _{i}$. Take as $F _{i}$ the extension of $W _{i}$ in
$\overline K _{i}$ generated by the $h _{i}$-th roots of $X _{i}$,
when $h _{i}$ runs across the set $\{\nu \in \mathbb{N}\colon \nu $
is not divisible by any $q _{i} \in \Pi _{i}\}$. The assertion that
$F _{i}$ is perfect and has the properties claimed by Lemma
\ref{lemm5.2} (b) is proved in the same way as the fact that the
fields defined by (4.1) satisfy the conditions of Proposition
\ref{prop4.1}. Since $\cup _{i=0} ^{m-1} \Pi _{i} = \mathbb{P}$ and
$\Pi _{i} \cap \Pi _{i'} = \emptyset $, $i \neq i'$, this implies
$F _{m-1}$ is quasifinite. Observing finally that ddim$(F _{i}) \le
{\rm ddim}(K _{i}) \le i + 1$, for $i = 1, \dots , m - 1$ (the latter
inequality follows from Greenberg's theorem), one completes the
proof.
\end{proof}
\par
\medskip
\begin{lemm}
\label{lemm5.3} Let $m$ be an integer $\ge 2$, and $\beta $  a
real number such that
\par\noindent
$\sqrt[m+1]{1 - (m + 1) ^{-2}} < \beta < 1$. Assume that $q \in
\overline {\mathbb{P}}$ and $\Pi _{i}\colon i = 0, \dots , m - 1$,
are nonempty subsets of $\mathbb{P}$ satisfying the following two
conditions:
\par
{\rm (c)} $q \notin \Pi _{m-1}$, $\cup _{i=0} ^{m-1} \Pi _{i} =
\mathbb{P}$, and $\Pi _{i'} \cap \Pi _{i''} = \emptyset $, for each
pair of indices $i' \neq i''$; in addition, $\Pi _{i}$ is finite
unless $i = m - 1$;
\par\vskip0.063truecm
{\rm (cc)} If $0 < i < m - 1$, then $\Pi _{i}$ consists of $(\beta
, i; q)$-representable numbers, and for each $q _{i} \in \Pi
_{i}$, we have $q _{i} = q _{i-1}' + q _{i-1}'' + q _{i-1}'''$,
for some $q _{i-1}'$, $q _{i-1}''$, $q _{i-1}''' \in \Pi _{i-1}$
with $q < q _{i} ^{\beta } < q _{i-1}' < q _{i-1}'' < q
_{i-1}'''$; also, $\Pi _{m-1}$ contains an element $q _{m-1}$
equal to $q _{m-2}' + q _{m-2}'' + q _{m-2}'''$, where $q
_{m-2}'$, $q _{m-2}''$ and $q _{m-2}'''$ are pairwise distinct
elements of $\Pi _{m-2}$, and $q < q _{m-1} ^{\beta } < {\rm
min}\{q _{m-2}', q _{m-2}'', q _{m-2}'''\}$.
\par\vskip0.063truecm
Then any quasifinite field $F _{m-1}$ singled out by Lemma
\ref{lemm5.2} admits a form $f _{m-1}$ of degree $q _{m-1}$ in
more than $q _{m-1} ^{m-1}$ variables, which does not possess a
nontrivial zero over $F _{m-1}$; in particular, {\rm ddim}$(F
_{m-1}) = m$.
\end{lemm}
\par
\medskip
\begin{proof}
Since, by Lemma \ref{lemm5.1} and the choice of $F _{m-1}$,
ddim$(F _{m-1}) \le m$, the latter part of our assertion follows
from the former one. Arguing by the method of proving Lemma
\ref{lemm4.2}, one obtains from conditions (c) and (cc) the
existence of an $F _{m-1}$-form $f _{m-1}$ of degree $q _{m-1}$ in
more than $q _{m-1} ^{\bar \lambda }$ variables, and without a
nontrivial zero over $F _{m-1}$, where $\bar \lambda = \sum _{i=0}
^{m-1} \beta ^{i}$. The assumptions on $\beta $ show that $\bar
\lambda > m - 1$, so Lemma \ref{lemm5.3} is proved.
\end{proof}
\par
\smallskip
It is now easy to complete the proof of Theorem \ref{theo2.2}.
Assume that $\ell < \infty $ and put $m = \ell - k + 1$. Applying
Lemma \ref{lemm4.3} (to the case where $\Pi = \mathbb{P}$) if $m =
1$, and Lemma \ref{lemm5.3} when $m \ge 2$, one obtains that there
is a quasifinite field $F _{m, q}$ with char$(F _{m,q}) = q$ and
ddim$(F _{m,q}) = m$, which admits an $F _{m, q}$-form $f _{m,q}$
of degree $q _{m-1} \in \mathbb{P} _{q}$ in at least $1 + q _{m-1}
^{m-1}$ variables, possessing no nontrivial zero over $F _{m-1}$.
Arguing by induction on $k - 1$, and using Lemma \ref{lemm3.2},
one can associate with $f _{m,q}$ an $E _{k, \ell } ^{\prime
}$-form $y _{m,q}$ of degree $q _{m-1}$ in more than $q _{m-1}
^{\ell - 1}$ variables, and without a nontrivial zero over $E _{k,
\ell } ^{\prime }$. Therefore, ddim$(E _{k, \ell } ^{\prime }) \ge
\ell $, and since ddim$(E _{k, \ell } ^{\prime }) \le {\rm ddim}(E
_{k, \ell }) = \ell $ (apply repeatedly Greenberg's theorem), we
have ddim$(E _{k, \ell } ^{\prime }) = \ell $, so Theorem
\ref{theo2.2} is proved.
\par
\medskip
As explained in Section~2, every perfect field $E$ satisfies
CD$(E) = {\rm cd}(\mathcal{G}_{E})$. Generally, this does not
apply to arbitrary fields (see (1.1) and (1.2)).  However, the
equality holds for interesting classes of imperfect fields, such
as the iterated Laurent series fields in $n$ variables, for any
fixed $n \in \mathbb{N}$, over quasifinite fields of type $C _{1}$
and nonzero characteristic. To demonstrate it, consider a field $E
_{n} = E _{0}((X _{1})) \dots ((X _{n}))$ of this kind over a
quasifinite $C _{1}$-field $E _{0}$ of characteristic $q$. Then
char$(E _{n}) = q$ and $[E _{n}\colon E _{n} ^{q}] = q ^{n}$,
which allows, similarly to the proof of Lemma \ref{lemm5.1}, to
deduce from Galois cohomology, combined with Kato-Milne cohomology
\cite[Theorem~3~(3)]{Ka}, Greenberg's theorem and the Arason-Baeza
theorem, that cd$_{q}(\mathcal{G}_{E _{n}}) = 1$, and the
dimensions cd$(\mathcal{G}_{E _{n}})$, CD$(E _{n})$, ddim$(E
_{n})$ and dim$_{q}(E _{n})$ (as well as $\overline {\rm CD}(E
_{n})$, see page \pageref{k45}) are equal to $n + 1$.
\par
\medskip
\begin{rema}
\label{rema5.4} Arguing as in the proof of Corollary
\ref{coro4.4}, one obtains that the class of quasifinite fields
with the property described by Lemma \ref{lemm5.3} is closed under
taking finite extensions. This implies in conjunction with
Corollary \ref{coro4.4} that the fields $E _{k, \ell ; q}$ singled
out by Theorem \ref{theo2.2} can be chosen so that their finite
extensions have Diophantine dimension $\ell $, for any $q \in
\overline {\mathbb{P}}$ and each pair $(k, \ell ) \neq (0, 0)$
admissible by Question 1. It is well-known that the open subgroups
of $\mathcal{G}_{E _{k, \ell ; q}}$, namely, the absolute Galois
groups of the considered fields have cohomological dimension $k$
(see \cite[Ch. I, Proposition~14]{S1}).
\end{rema}
\par
\medskip
To conclude with, let us note that our proof of Theorem
\ref{theo2.1} strongly depends on the fact that given a
quasifinite field $F$, the Sylow pro-$p$-subgroup of $\mathcal{G}
_{F}$ is isomorphic to $\mathbb{Z} _{p}$, for every $p \in
\mathbb{P}$. The existence of perfect fields $E _{p}$, $p \in
\mathbb{P} _{2}$, with dim$(E _{p}) \le 1$, ddim$(E _{p}) \ge 2$
and $\mathcal{G}_{E _{p}}$ a pro-$p$-group with infinitely many
open subgroups of index $p$, solving \cite[Problem~2]{KK}, has
been established in \cite{CTM} and \cite{C-T}, by a different
method. Both results (as well as Merkur'ev's counter-example to
the Kato-Kuzumaki $C _{2} ^{0}$ conjecture, see \cite{KK} and
\cite[Theorem~4 and Proposition~5]{Me}) leave open the following
question.
\par
\medskip
{\bf Question 2.} Find whether there exists a perfect field $E$
such that dim$(E) = \infty $ and some of the following two
conditions holds:
\par
(a) dim$(E) \le 1$ and the Sylow pro-$p$-subgroups of
$\mathcal{G}_{E}$ are not isomorphic to $\mathbb{Z} _{p}$, for any
$p \in \mathbb{P}$ (with, possibly, finitely many exceptions);
\par
(b) $E$ is the fixed field of a Sylow pro-$p$-subgroup $G _{p}$ of
$\mathcal{G}_{\mathbb{Q}_{1}}$, for some $p \in \mathbb{P}$, where
$\mathbb{Q} _{1}$ is a totally imaginary number field.
\par
\medskip
As noted at the end of Remark \ref{rema2.3}, nothing seems to be
known about the Diophantine dimensions of the fixed fields of
Sylow's subgroups of $\mathcal{G}_{\mathbb{Q} _{1}}$, for any
totally imaginary number field $\mathbb{Q} _{1}$. Against this
background, Question 2 (b) makes interest in its own right. A
negative answer to it, for some $p \in \mathbb{P}$, would show
that abrd$_{p}(\mathbb{Q} ^{\prime }) < \infty $, for all
finitely-generated extensions $\mathbb{Q} ^{\prime }/\mathbb{Q}$.
This follows from the Lang-Nagata-Tsen theorem, the main results
of \cite{Mat}, and \cite[Corollary~4.5.11]{GiSz}.
\par
\vskip0.32truecm\noindent \emph{Acknowledgements.} The authors
wish to thank the referee for several useful suggestions and
comments, particularly, for drawing their attention to the
Arason-Baeza theorem and its relation with Question 1, which were
used for improving the quality of the present paper. The research
of the first-named author has partially been supported by the
Bulgarian National Science Fund under Grant KP-06 N 32/1 of
07.12.2019; the research of the second-named author has partially
been supported by Sofia University contract 80-10-99/2023.

\end{document}